\documentclass[journal, twoside]{IEEEtran}

\usepackage{import}

\usepackage[utf8]{luainputenc}

\usepackage[T1]{fontenc}

\usepackage{graphicx}
\usepackage[dvipsnames]{xcolor}

\usepackage[english]{babel}
\addto\captionsenglish{}
\addto\captionsenglish{}
\usepackage{csquotes}

\usepackage{newtxtext}
\usepackage{amsthm}
\usepackage[slantedGreek]{newtxmath}
\usepackage[OMLmathsfit]{isomath}
\DeclareMathAlphabet{\mathbfsf}{\encodingdefault}{\sfdefault}{bx}{n}
\usepackage{bm}
\usepackage{envmath}
\usepackage{mathtools}
\usepackage{commath}
\usepackage{siunitx}

\usepackage{subcaption}
\makeatletter
\providecommand{\p@subfigure}{}

\renewcommand{\p@subfigure}{\thefigure}

\makeatother

\usepackage{booktabs}
\usepackage{footmisc}  

\usepackage{url}

\theoremstyle{definition}

\theoremstyle{plain}

\theoremstyle{remark}

\usepackage{lineno}
\modulolinenumbers[5]
\usepackage{todonotes}
\usepackage{umoline}

\usepackage{pgfplots}
\usepackage{pgfplotstable}
\pgfplotsset{compat=newest}
\pgfplotsset{plot coordinates/math parser=false}
\newlength\figureheight
\newlength\figurewidth
\pgfplotsset{every axis plot/.append style={line width=1.5pt},
    legend style={font=\footnotesize, 
        text height=1.0ex,
        draw=black,
        fill=white,
        legend cell align=left}}

\usepackage[hidelinks]{hyperref} 
\usepackage[english]{cleveref}

\Crefname{defn}{definition}{definitions}
\Crefname{defn}{Definition}{Definitions}

\Crefname{asm}{assumption}{assumptions}
\Crefname{asm}{Assumption}{Assumptions}

\crefname{lem}{lemma}{lemmas} 
\Crefname{lem}{Lemma}{Lemmas}

\crefname{prop}{proposition}{propositions} 
\Crefname{prop}{Proposition}{Propositions}

\crefname{thm}{theorem}{theorms} 
\Crefname{thm}{Theorem}{Theorms}

\crefname{cor}{corollary}{corollaries}
\Crefname{cor}{Corollary}{Corollaries}

\Crefname{figure}{Fig.}{Figs.}

\AtBeginEnvironment{appendices}{\crefalias{section}{appendix}}

\usepackage{upgreek}

\newcounter{subequation}
\newlength\mtabskip\mtabskip=-1.25cm

\def\mtabLong{long}

\newcommand{\mr}{\mathrm}

\newcommand{\mc}{\mathcal}

\newcommand{\veg}[1]{\bm{#1}}     
\newcommand{\mat}[1]{%
  \ifx#1U%
    \mathsfbfit{U}\mkern1mu%
  \else\ifx#1V%
    \mathsfbfit{V}\mkern1mu%
  \else\ifx#1\Theta%
    \mathsfbfit{\Theta}\mkern0.5mu%
  \else%
    \mathsfbfit{#1}%
  \fi\fi\fi%
}
\renewcommand{\vec}[1]{\mathsfbfit{#1}} 
\newcommand{\wmat}[1]{\widetilde{\mathsfbfit{#1}}} 
\newcommand{\op}[1]{\mathcal{#1}} 
\newcommand{\vecop}[1]{\bm{\mathcal{#1}}} 

\newcommand{\n}{\hat{\bm{n}}}

\newcommand{\dd}{\mathrm{d}}  

\newcommand{\jm}{\mathrm{j}}  

\newcommand{\e}{\mathrm{e}}

\newcommand{\normF}[1]{\norm{#1}_\mathrm{F}}

\DeclareMathOperator{\dist}{dist}
\DeclareMathOperator{\diam}{diam}

\DeclareMathOperator*{\argmax}{arg\,max}

\newcommand{\duality}[2]{\langle #1, #2\rangle_\Gamma}

\newcommand\restr[2]{{
        \left.\kern-\nulldelimiterspace 
        #1 
        \vphantom{|} 
        \right|_{#2} 
}}

\newcommand\rst[3]{{
        \left.\kern-\nulldelimiterspace 
        #1 
        \vphantom{|} 
        \right|_{#2}^{#3} 
}}

\usepackage{acro}
\usepackage{pdfcomment}
\acsetup{pdfcomments/use=true}

\DeclareAcronym{DG}
{
    short = DG ,
    long = discontinuous Galerkin
}

\DeclareAcronym{ACA}
{
    short = ACA ,
    long = adaptive cross approximation
}

\DeclareAcronym{EFIE}
{
    short =  EFIE ,
    long = electric field integral equation
}

\DeclareAcronym{MFIE}
{
    short =  MFIE ,
    long = magnetic field integral equation
}

\DeclareAcronym{CFIE}
{
    short =  CFIE ,
    long = combined field integral equation
}

\DeclareAcronym{MUIE}
{
    short =  MUIE ,
    long = Müller integral equation
}

\DeclareAcronym{PMCHWT}
{
    short =  PMCHWT ,
    long = Poggio-Miller-Chang-Harrington-Wu-Tsai integral equation
}

\DeclareAcronym{MLFMA}
{
    short =  MLFMA ,
    long = multilevel fast multipole algorithm
}

\DeclareAcronym{AIM}
{
    short =  AIM ,
    long = adaptive integral method
}

\DeclareAcronym{BEM}
{
    short =  BEM ,
    long = boundary element method
}

\DeclareAcronym{SPD}
{
    short =  SPD ,
    long = {symmetric, positive definite}
}

\DeclareAcronym{SPSD}
{
    short =  SPD ,
    long = {symmetric, positive semi-definite}
}

\DeclareAcronym{PEC}
{
    short =  PEC ,
    long = perfectly electrically conducting
}

\DeclareAcronym{RWG}
{
    short = RWG ,
    long = Rao-Wilton-Glisson
} 

\DeclareAcronym{BC}
{
    short = BC ,
    long = Buffa-Christiansen
}

\DeclareAcronym{SVD}
{
    short = SVD ,
    long = singular value decomposition
}

\DeclareAcronym{CG}
{
    short = CG ,
    long = conjugate gradient
} 

\DeclareAcronym{PCG}
{
    short = PCG ,
    long = preconditioned conjugate gradient
} 

\DeclareAcronym{CGS}
{
    short = CGS ,
    long = conjugate gradient squared
}

\DeclareAcronym{CMP}
{
    short = CMP ,
    long = Calderón multiplicative preconditioner
} 

\DeclareAcronym{RFCMP}
{
    short = RF-CMP ,
    long = refinement-free Calderón multiplicative preconditioner
} 

\DeclareAcronym{HPD}
{
    short = HPD ,
    long = {Hermitian, positive definite}
} 

\DeclareAcronym{RHS}
{
    short = RHS ,
    long = right-hand side
}

\DeclareAcronym{LSE}
{
    short = LSE ,
    long = linear system of equations
}

\DeclareAcronym{AMG}
{
    short = AMG ,
    long = algebraic multigrid
}

\DeclareAcronym{PW}
{
    short = PW ,
    long = plane wave
}

\DeclareAcronym{GMRES}
{
    short = GMRES ,
    long = generalized minimum residual
}

\DeclareAcronym{IDR}
{
    short = IDR ,
    long = induced dimension reduction
}

\DeclareAcronym{BICGstab}
{
    short = BiCGstab ,
    long = stabilized bi-conjugate gradient
}

\DeclareAcronym{FF}
{
    short = FF ,
    long = far field
}

\DeclareAcronym{NF}
{
    short = NF ,
    long = near field
}

\DeclareAcronym{SC}
{
    short = SCC ,
    long = standard convergence criterion ,
    pdfcomment = standard convergence criterion
}

\DeclareAcronym{RC}
{
    short = RSCC ,
    long = random-sampling convergence criterion ,
    pdfcomment = random-sampling convergence criterion
}

\DeclareAcronym{CC}
{
    short = CCC ,
    long = combined convergence criterion ,
    pdfcomment = combined convergence criterion
}

\DeclareAcronym{HCA}
{
    short = HCA,
    long = hybrid cross approximation
}

\DeclareAcronym{FMM}
{
    short = FMM ,
    long = fast multipole method  ,
    pdfcomment = fast multipole method 
}

\DeclareAcronym{HMatrix}
{
    short = $\mathcal{H}$-matrix ,
    long = hierarchical matrix ,
    pdfcomment = hierarchical matrix 
}

\DeclareAcronym{MLMDA}
{
    short = MLMDA ,
    long = multilevel matrix decomposition algorithm ,
    pdfcomment = multilevel matrix decomposition algorithm
}

\DeclareAcronym{NCA}
{
    short = NCA ,
    long = nested cross approximation ,
    pdfcomment = nested cross approximation
}

\DeclareAcronym{NPSA}
{
    short = NPSA ,
    long = nested pseudo skeleton approximation ,
    pdfcomment = nested pseudo skeleton approximation
}

\DeclareAcronym{IACA}
{
    short = IACA ,
    long = incomplete adaptive cross approximation ,
    pdfcomment = incomplete adaptive cross approximation ,
}

\DeclareAcronym{SIE}
{
    short = SIE ,
    long = surface integral equation,
    pdfcomment = surface integral equation ,
}

\DeclareAcronym{MVP}
{
    short = MVP ,
    long = matrix--vector product,
    pdfcomment = matrix--vector product
}

\DeclareAcronym{MoM}
{
    short = MoM ,
    long = method of moments,
    pdfcomment = method of moments
}

\DeclareAcronym{NESA}
{
    short = NESA ,
    long = nested equivalent source approximation ,
    pdfcomment = nested equivalent source approximation 
}

\DeclareAcronym{WNESA}
{
    short = WNESA ,
    long = wideband nested equivalent source approximation ,
    pdfcomment = wideband nested equivalent source approximation 
}

\DeclareAcronym{RCS}
{
    short = RCS ,
    long = radar cross section ,
    pdfcomment = radar cross section 
}

\DeclareAcronym{DH2}
{
    short = {$\mc D\mc H^2$-matrix} ,
    short-plural-form = {$\mc D\mc H^2$-matrices} ,
    long = {directional $\mc H^2$-matrix} ,
    long-plural-form = {directional $\mc H^2$-matrices} ,
    pdfstring = {directional H2-matrix} ,
} 

\newcolumntype {n}{c}
\newcolumntype {N}{>{\small}c}
\newcolumntype {L}{>{\small}l}
\newcolumntype {F}{>{\footnotesize}c}
\newcolumntype {v}[1]{>{\raggedright \hspace {0pt}} p {#1}}
\newcolumntype {V}[1]{>{\small \raggedright \hspace {0pt}} p {#1}}
\newcolumntype{d}[1]{>{\DC@{.}{.}{#1}}c<{\DC@end}}

\newcolumntype{R}[1]{%
    >{\begin{turn}{90}\begin{minipage}{#1}\small\raggedright\hspace{0pt}}l%
            <{\end{minipage}\end{turn}}%
}

\pgfplotsset{colormap={hawaii}{
		rgb = (0.550541, 0.006842, 0.451980);
		rgb = (0.551494, 0.015367, 0.447972);
		rgb = (0.552426, 0.023795, 0.443998);
		rgb = (0.553328, 0.032329, 0.440021);
		rgb = (0.554227, 0.041170, 0.436063);
		rgb = (0.555098, 0.049286, 0.432125);
		rgb = (0.555948, 0.056667, 0.428188);
		rgb = (0.556797, 0.063525, 0.424272);
		rgb = (0.557619, 0.069970, 0.420377);
		rgb = (0.558415, 0.076028, 0.416509);
		rgb = (0.559210, 0.081936, 0.412663);
		rgb = (0.559991, 0.087507, 0.408823);
		rgb = (0.560746, 0.092811, 0.405012);
		rgb = (0.561495, 0.098081, 0.401237);
		rgb = (0.562235, 0.103128, 0.397471);
		rgb = (0.562954, 0.108005, 0.393736);
		rgb = (0.563663, 0.112872, 0.390025);
		rgb = (0.564355, 0.117530, 0.386344);
		rgb = (0.565032, 0.122122, 0.382698);
		rgb = (0.565709, 0.126681, 0.379074);
		rgb = (0.566380, 0.131171, 0.375474);
		rgb = (0.567037, 0.135542, 0.371905);
		rgb = (0.567679, 0.139872, 0.368378);
		rgb = (0.568312, 0.144198, 0.364861);
		rgb = (0.568939, 0.148416, 0.361384);
		rgb = (0.569559, 0.152618, 0.357942);
		rgb = (0.570171, 0.156806, 0.354519);
		rgb = (0.570777, 0.160934, 0.351127);
		rgb = (0.571377, 0.165008, 0.347764);
		rgb = (0.571972, 0.169120, 0.344417);
		rgb = (0.572562, 0.173131, 0.341120);
		rgb = (0.573142, 0.177166, 0.337836);
		rgb = (0.573711, 0.181138, 0.334602);
		rgb = (0.574276, 0.185151, 0.331356);
		rgb = (0.574840, 0.189095, 0.328170);
		rgb = (0.575406, 0.193035, 0.324992);
		rgb = (0.575967, 0.196978, 0.321854);
		rgb = (0.576518, 0.200854, 0.318740);
		rgb = (0.577060, 0.204783, 0.315654);
		rgb = (0.577596, 0.208664, 0.312565);
		rgb = (0.578135, 0.212545, 0.309542);
		rgb = (0.578676, 0.216431, 0.306516);
		rgb = (0.579214, 0.220287, 0.303496);
		rgb = (0.579746, 0.224106, 0.300518);
		rgb = (0.580271, 0.227977, 0.297566);
		rgb = (0.580793, 0.231817, 0.294618);
		rgb = (0.581315, 0.235646, 0.291715);
		rgb = (0.581835, 0.239463, 0.288810);
		rgb = (0.582353, 0.243268, 0.285910);
		rgb = (0.582870, 0.247097, 0.283066);
		rgb = (0.583386, 0.250916, 0.280201);
		rgb = (0.583901, 0.254739, 0.277381);
		rgb = (0.584416, 0.258531, 0.274552);
		rgb = (0.584931, 0.262342, 0.271740);
		rgb = (0.585443, 0.266156, 0.268980);
		rgb = (0.585951, 0.269966, 0.266198);
		rgb = (0.586456, 0.273771, 0.263439);
		rgb = (0.586961, 0.277575, 0.260676);
		rgb = (0.587466, 0.281374, 0.257925);
		rgb = (0.587972, 0.285180, 0.255221);
		rgb = (0.588478, 0.289013, 0.252494);
		rgb = (0.588984, 0.292818, 0.249767);
		rgb = (0.589491, 0.296652, 0.247081);
		rgb = (0.589999, 0.300465, 0.244376);
		rgb = (0.590507, 0.304300, 0.241716);
		rgb = (0.591016, 0.308135, 0.239031);
		rgb = (0.591526, 0.311969, 0.236379);
		rgb = (0.592038, 0.315846, 0.233692);
		rgb = (0.592548, 0.319698, 0.231058);
		rgb = (0.593055, 0.323559, 0.228420);
		rgb = (0.593562, 0.327429, 0.225773);
		rgb = (0.594071, 0.331309, 0.223134);
		rgb = (0.594583, 0.335229, 0.220510);
		rgb = (0.595095, 0.339131, 0.217865);
		rgb = (0.595609, 0.343048, 0.215226);
		rgb = (0.596126, 0.346976, 0.212613);
		rgb = (0.596645, 0.350921, 0.209994);
		rgb = (0.597164, 0.354880, 0.207388);
		rgb = (0.597680, 0.358830, 0.204776);
		rgb = (0.598196, 0.362821, 0.202147);
		rgb = (0.598721, 0.366829, 0.199533);
		rgb = (0.599248, 0.370837, 0.196964);
		rgb = (0.599771, 0.374879, 0.194370);
		rgb = (0.600294, 0.378931, 0.191738);
		rgb = (0.600819, 0.383009, 0.189149);
		rgb = (0.601346, 0.387090, 0.186548);
		rgb = (0.601874, 0.391215, 0.183949);
		rgb = (0.602403, 0.395345, 0.181345);
		rgb = (0.602933, 0.399486, 0.178782);
		rgb = (0.603464, 0.403678, 0.176158);
		rgb = (0.603995, 0.407873, 0.173594);
		rgb = (0.604521, 0.412102, 0.171015);
		rgb = (0.605043, 0.416348, 0.168436);
		rgb = (0.605562, 0.420618, 0.165848);
		rgb = (0.606084, 0.424928, 0.163317);
		rgb = (0.606609, 0.429252, 0.160731);
		rgb = (0.607129, 0.433600, 0.158195);
		rgb = (0.607639, 0.437998, 0.155649);
		rgb = (0.608144, 0.442412, 0.153086);
		rgb = (0.608644, 0.446848, 0.150582);
		rgb = (0.609134, 0.451324, 0.148071);
		rgb = (0.609610, 0.455826, 0.145615);
		rgb = (0.610079, 0.460356, 0.143119);
		rgb = (0.610542, 0.464933, 0.140685);
		rgb = (0.610991, 0.469544, 0.138267);
		rgb = (0.611421, 0.474170, 0.135829);
		rgb = (0.611833, 0.478839, 0.133514);
		rgb = (0.612226, 0.483539, 0.131212);
		rgb = (0.612600, 0.488287, 0.128920);
		rgb = (0.612950, 0.493049, 0.126718);
		rgb = (0.613275, 0.497875, 0.124574);
		rgb = (0.613572, 0.502705, 0.122487);
		rgb = (0.613837, 0.507592, 0.120512);
		rgb = (0.614069, 0.512502, 0.118669);
		rgb = (0.614264, 0.517459, 0.116848);
		rgb = (0.614418, 0.522434, 0.115160);
		rgb = (0.614530, 0.527456, 0.113657);
		rgb = (0.614594, 0.532510, 0.112266);
		rgb = (0.614607, 0.537595, 0.111032);
		rgb = (0.614566, 0.542708, 0.109999);
		rgb = (0.614468, 0.547849, 0.109114);
		rgb = (0.614308, 0.553016, 0.108421);
		rgb = (0.614082, 0.558212, 0.108010);
		rgb = (0.613787, 0.563446, 0.107850);
		rgb = (0.613419, 0.568682, 0.107943);
		rgb = (0.612974, 0.573946, 0.108312);
		rgb = (0.612449, 0.579232, 0.109026);
		rgb = (0.611842, 0.584522, 0.110040);
		rgb = (0.611148, 0.589820, 0.111320);
		rgb = (0.610353, 0.595132, 0.112963);
		rgb = (0.609471, 0.600443, 0.114856);
		rgb = (0.608494, 0.605748, 0.117169);
		rgb = (0.607411, 0.611060, 0.119811);
		rgb = (0.606215, 0.616350, 0.122763);
		rgb = (0.604930, 0.621618, 0.126124);
		rgb = (0.603536, 0.626876, 0.129757);
		rgb = (0.602026, 0.632107, 0.133692);
		rgb = (0.600413, 0.637306, 0.137967);
		rgb = (0.598689, 0.642469, 0.142496);
		rgb = (0.596862, 0.647588, 0.147334);
		rgb = (0.594916, 0.652662, 0.152416);
		rgb = (0.592872, 0.657697, 0.157790);
		rgb = (0.590707, 0.662667, 0.163419);
		rgb = (0.588441, 0.667579, 0.169258);
		rgb = (0.586085, 0.672429, 0.175280);
		rgb = (0.583613, 0.677213, 0.181507);
		rgb = (0.581049, 0.681916, 0.187985);
		rgb = (0.578388, 0.686560, 0.194586);
		rgb = (0.575646, 0.691121, 0.201310);
		rgb = (0.572809, 0.695614, 0.208243);
		rgb = (0.569878, 0.700018, 0.215285);
		rgb = (0.566888, 0.704346, 0.222470);
		rgb = (0.563814, 0.708597, 0.229738);
		rgb = (0.560662, 0.712753, 0.237171);
		rgb = (0.557458, 0.716845, 0.244622);
		rgb = (0.554182, 0.720839, 0.252219);
		rgb = (0.550853, 0.724766, 0.259874);
		rgb = (0.547470, 0.728605, 0.267574);
		rgb = (0.544043, 0.732376, 0.275394);
		rgb = (0.540571, 0.736058, 0.283238);
		rgb = (0.537067, 0.739685, 0.291141);
		rgb = (0.533507, 0.743228, 0.299094);
		rgb = (0.529936, 0.746702, 0.307079);
		rgb = (0.526333, 0.750112, 0.315113);
		rgb = (0.522696, 0.753461, 0.323192);
		rgb = (0.519049, 0.756752, 0.331281);
		rgb = (0.515367, 0.759983, 0.339437);
		rgb = (0.511681, 0.763162, 0.347595);
		rgb = (0.507990, 0.766293, 0.355785);
		rgb = (0.504280, 0.769372, 0.363984);
		rgb = (0.500550, 0.772410, 0.372217);
		rgb = (0.496820, 0.775405, 0.380485);
		rgb = (0.493085, 0.778365, 0.388763);
		rgb = (0.489350, 0.781287, 0.397049);
		rgb = (0.485614, 0.784180, 0.405376);
		rgb = (0.481884, 0.787038, 0.413711);
		rgb = (0.478142, 0.789866, 0.422057);
		rgb = (0.474411, 0.792674, 0.430440);
		rgb = (0.470680, 0.795455, 0.438824);
		rgb = (0.466955, 0.798219, 0.447235);
		rgb = (0.463220, 0.800964, 0.455667);
		rgb = (0.459518, 0.803693, 0.464121);
		rgb = (0.455810, 0.806409, 0.472577);
		rgb = (0.452124, 0.809110, 0.481054);
		rgb = (0.448436, 0.811796, 0.489555);
		rgb = (0.444772, 0.814472, 0.498091);
		rgb = (0.441108, 0.817144, 0.506616);
		rgb = (0.437487, 0.819803, 0.515175);
		rgb = (0.433858, 0.822465, 0.523755);
		rgb = (0.430280, 0.825110, 0.532352);
		rgb = (0.426720, 0.827756, 0.540960);
		rgb = (0.423186, 0.830401, 0.549598);
		rgb = (0.419708, 0.833036, 0.558241);
		rgb = (0.416257, 0.835673, 0.566923);
		rgb = (0.412868, 0.838305, 0.575612);
		rgb = (0.409520, 0.840937, 0.584314);
		rgb = (0.406245, 0.843562, 0.593044);
		rgb = (0.403035, 0.846190, 0.601780);
		rgb = (0.399905, 0.848819, 0.610541);
		rgb = (0.396872, 0.851439, 0.619320);
		rgb = (0.393950, 0.854061, 0.628104);
		rgb = (0.391152, 0.856683, 0.636905);
		rgb = (0.388472, 0.859301, 0.645709);
		rgb = (0.385935, 0.861918, 0.654530);
		rgb = (0.383585, 0.864526, 0.663367);
		rgb = (0.381407, 0.867128, 0.672196);
		rgb = (0.379424, 0.869728, 0.681023);
		rgb = (0.377672, 0.872325, 0.689863);
		rgb = (0.376170, 0.874907, 0.698686);
		rgb = (0.374923, 0.877482, 0.707507);
		rgb = (0.373981, 0.880045, 0.716318);
		rgb = (0.373340, 0.882596, 0.725106);
		rgb = (0.373043, 0.885136, 0.733865);
		rgb = (0.373112, 0.887654, 0.742601);
		rgb = (0.373570, 0.890156, 0.751300);
		rgb = (0.374439, 0.892639, 0.759946);
		rgb = (0.375723, 0.895095, 0.768546);
		rgb = (0.377467, 0.897524, 0.777098);
		rgb = (0.379671, 0.899923, 0.785572);
		rgb = (0.382352, 0.902288, 0.793974);
		rgb = (0.385527, 0.904619, 0.802283);
		rgb = (0.389213, 0.906913, 0.810503);
		rgb = (0.393385, 0.909161, 0.818619);
		rgb = (0.398074, 0.911369, 0.826627);
		rgb = (0.403255, 0.913528, 0.834507);
		rgb = (0.408926, 0.915628, 0.842255);
		rgb = (0.415083, 0.917688, 0.849859);
		rgb = (0.421704, 0.919678, 0.857309);
		rgb = (0.428791, 0.921615, 0.864606);
		rgb = (0.436305, 0.923489, 0.871734);
		rgb = (0.444231, 0.925293, 0.878682);
		rgb = (0.452541, 0.927032, 0.885454);
		rgb = (0.461203, 0.928705, 0.892037);
		rgb = (0.470211, 0.930311, 0.898424);
		rgb = (0.479521, 0.931839, 0.904620);
		rgb = (0.489103, 0.933297, 0.910617);
		rgb = (0.498950, 0.934685, 0.916408);
		rgb = (0.509019, 0.936004, 0.922005);
		rgb = (0.519281, 0.937246, 0.927394);
		rgb = (0.529715, 0.938416, 0.932588);
		rgb = (0.540292, 0.939517, 0.937592);
		rgb = (0.550997, 0.940549, 0.942401);
		rgb = (0.561804, 0.941509, 0.947020);
		rgb = (0.572686, 0.942411, 0.951459);
		rgb = (0.583621, 0.943243, 0.955728);
		rgb = (0.594606, 0.944015, 0.959825);
		rgb = (0.605610, 0.944731, 0.963765);
		rgb = (0.616637, 0.945388, 0.967563);
		rgb = (0.627648, 0.945989, 0.971214);
		rgb = (0.638645, 0.946543, 0.974739);
		rgb = (0.649620, 0.947052, 0.978146);
		rgb = (0.660548, 0.947515, 0.981449);
		rgb = (0.671439, 0.947934, 0.984653);
		rgb = (0.682276, 0.948316, 0.987765);
		rgb = (0.693064, 0.948662, 0.990803);
		rgb = (0.703779, 0.948977, 0.993775);
}}

\NewDocumentCommand{\TT}{o}{
    \IfNoValueTF {#1} {%
        \vecop T
    }
    {
        \vecop T_{#1}
    }
}

\NewDocumentCommand{\TA}{o}{
    \IfNoValueTF {#1} {%
        \vecop T_{\kern-2pt\mr{A}}
    }
    {
        \vecop T_{\kern-2pt\mr{A},#1}
    }
}

\NewDocumentCommand{\V}{o}{
    \IfNoValueTF {#1} {%
        \op V
    }
    {
        \op V_{#1}
    }
}

\NewDocumentCommand{\W}{o}{
    \IfNoValueTF {#1} {%
        \op W
    }
    {
        \op W_{#1}
    }
}

\NewDocumentCommand{\TPhi}{o}{
    \IfNoValueTF {#1} {%
        \vecop T_{\kern-2pt\Phiup}
    }
    {
        \vecop T_{\kern-2pt\Phiup,#1}
    }
}

\NewDocumentCommand{\matTA}{o}{
    \IfNoValueTF {#1} {%
        \mat T_\mr{A}   
        }
    {
        \mat T_{\mr{A},#1}
    }
}

\NewDocumentCommand{\matTPhi}{o}{
    \IfNoValueTF {#1} {%
        \mat T_\Phiup
        }
    {
        \mat T_{\Phiup,#1}
    }
}

\NewDocumentCommand{\matV}{o}{
    \IfNoValueTF {#1} {%
        \mat V   
        }
    {
        \mat V_{#1}
    }
}

\newcommand{\fk}{f\kern-2pt} 
\usepackage[style=ieee, backend=biber, isbn=false, maxbibnames=9]{biblatex}

\usepackage{nicefrac}
\usepackage{yfonts}
\usepackage{comment}
\usepackage{algpseudocode}
\usepackage{algorithm}
\graphicspath{{figures/}}
\DeclareUnicodeCharacter{2212}{\ensuremath{-}}

\newcommand{\edits}[1]{}

\begin{document}
  \title{Wideband Directional $\mc H^2$-Matrix Compression for the Electric Field Integral Equation with Geometry-Adaptive Cluster Trees}
  
  %
  
  \author{Joshua M.\ Tetzner,~\IEEEmembership{Graduate Student Member,~IEEE,}
    and~Simon B.\ Adrian,~\IEEEmembership{Senior Member,~IEEE}
    \thanks{Funded by the Bundesministerium für Forschung, Technologie und Raumfahrt (BMFTR, Federal Ministry of Research, Technology and Space) under Grant 05H24HRA. \emph{(Corresponding author: Simon B. Adrian.)}}
    \thanks{J.\ M.\ Tetzner and S.\ B.\ Adrian are with the Fakultät für Informatik und Elektrotechnik, Universität Rostock, 18059 Rostock, Germany.}
  }

\maketitle

\begin{abstract}
  We present an efficient wideband construction of \aclp{DH2} for the \acl{EFIE} that, in contrast to existing constructions, supports not only box trees but also geometry-adaptive cluster trees. 
  To accommodate geometry-adaptive cluster trees, we determine the number of directions from the electrical size of each cluster (instead of the level of a box tree), construct the directions using Spherical-Fibonacci points, and establish a hierarchy between the direction sets of a cluster and its children through an angular nearest-neighbor mapping.
  We construct the nested directional representation using the \acl{IACA}, for which we introduce a robustified tree-mimicry pivoting strategy that prevents premature convergence for block-structured matrices arising from certain geometries and meshes.
  Numerical results demonstrate that the proposed approach achieves the desired accuracy, requires no more storage than the octree-based construction and substantially less when the geometry or discretization is poorly matched to octree clustering, and exhibits the expected $\mc O(N\log N)$ scaling for high-frequency problems.
\end{abstract}

\acresetall


%
\IEEEpeerreviewmaketitle

\section{Introduction}

\IEEEPARstart{E}{lectromagnetic} scattering and radiation problems are commonly formulated by surface integral equations, such as the \ac{EFIE}, \ac{MFIE}, or combined field integral equation, and discretized by the \ac{MoM} \cite{harrington_FieldComputationMoment_1993}.
The resulting system matrix is dense, leading to $\mc O(N^2)$ complexity for the setup time, storage cost, and the cost of a single \ac{MVP}, where $N$ denotes the number of unknowns.
For practically relevant problem sizes, fast methods are therefore indispensable.

Established methods include expansion-based approaches such as the \ac{MLFMA} \cite{coifman_FastMultipoleMethod_1993, lu_MultilevelAlgorithmSolving_1994, song_MultilevelFastMultipole_1997}, FFT--accelerated techniques such as the \ac{AIM} \cite{bleszynski_AIMAdaptiveIntegral_1996}, butterfly-type algorithms such as the \ac{MLMDA} \cite{michielssen_MultilevelMatrixDecomposition_1996}, directional multilevel algorithms for oscillatory kernels \cite{engquist_FastDirectionalMultilevel_2007}, wideband kernel-independent methods such as the \ac{WNESA}~\cite{li_WidebandFastKernelIndependent_2015}, and the mixed-form nested approximation~\cite{li_MixedformNestedApproximation_2018}.
A complementary class of fast methods is based on algebraic matrix compression, which constructs compressed representations directly from selected entries of the system matrix.
The \ac{ACA} is a widely used representative and approximates suitable matrix blocks by low-rank factorizations with ranks close to the optimal numerical rank.
The efficiency of such low-rank approximations, however, depends on the electrical size of the problem.

For electrically small problems, interactions of sufficiently separated clusters admit low-rank approximations with ranks bounded  independently of the wavenumber; such interactions are called admissible.
$\mc H$-matrices compress each admissible interaction independently and reach $\mc O(N\log N)$ complexity for storage, assembly, and a single \ac{MVP} \cite{zhao_AdaptiveCrossApproximation_2005}.
Nested low-rank representations of these blocks turn the $\mc H$-matrix into an $\mc H^2$-matrix and reduce the storage and \ac{MVP} cost to $\mc O(N)$.
The \ac{NCA} constructs such nested representations algebraically~\cite{bebendorf_ConstructingNestedBases_2012}.
Incorporating geometrical information into the \ac{NCA} and replacing its \ac{ACA} with the recently introduced \ac{IACA} makes this construction entry-efficient and reduces its assembly complexity to the near-optimal $\mc O(N)$ \cite{tetzner_IncompleteAdaptiveCross_2026}.

For electrically large problems, separation alone no longer bounds the rank: the ranks of the corresponding low-rank approximations grow with the electrical size of the associated clusters, so constructions based only on geometrical separation, such as the \ac{NCA}, lose their favorable rank behavior for electrically large clusters \cite{bebendorf_ConstructingNestedBases_2012}.
Directional methods restore bounded ranks by combining a wavenumber-dependent admissibility criterion with a subdivision of each admissible region into cone-like directional subsets \cite{engquist_FastDirectionalMultilevel_2007,bebendorf_WidebandNestedCross_2015}.
For sufficiently narrow directional subsets, the analysis of the free-space Green's function shows that its oscillatory behavior can be represented by a plane wave along the associated direction, while the remaining dependence is asymptotically smooth \cite{bebendorf_WidebandNestedCross_2015,borm_Directional$mathcalH^2$matrix_2017}.
Consequently, the corresponding interactions again admit low-rank approximations with ranks bounded independently of the wavenumber, which can be organized in a \ac{DH2} representation \cite{bebendorf_WidebandNestedCross_2015, borm_Directional$mathcalH^2$matrix_2017, noor_ApplicationDH2MatrixFramework_2026}.

Within the hierarchy of an electrically large problem, however, the electrical size of the clusters decreases toward the leaves, so electrically large and small clusters generally occur within the same cluster tree.
A wideband \ac{DH2} construction therefore uses directional nested low-rank approximations for electrically large clusters and standard nested low-rank approximations for electrically small clusters.
The electrical size of each cluster, therefore, determines the transition between the two representations \cite{bebendorf_WidebandNestedCross_2015, borm_Directional$mathcalH^2$matrix_2017}.

Existing \ac{DH2} constructions \cite{bebendorf_WidebandNestedCross_2015,borm_Directional$mathcalH^2$matrix_2017} prescribe their directions by a recursive face refinement of a cube and tie the resulting sets of directions to the levels of the cluster tree.
This fits regularly refined box trees, where a cluster's diameter is determined by its level.
Such trees can be highly effective when their fixed box subdivision partitions the unknowns approximately uniformly, as for a cube aligned with an octree.
For geometries or discretizations that are poorly matched to this subdivision, clusters on the same level may contain substantially different numbers of unknowns because some boxes intersect only small parts of the surface.
Geometry-adaptive cluster trees may follow the distribution of the unknowns more closely and thereby improve the admissible partition and the resulting compression.

However, geometry-adaptive trees may place clusters of substantially different diameter on the same level, so the required directional resolution can no longer be prescribed level-wise but must follow the electrical size of each cluster. 

In this article, we present a geometry-adaptive wideband construction of \acp{DH2} for the \ac{EFIE} that handles electrically large and small clusters in a single bottom-up traversal.
To this end, we combine (i) a geometry-adaptive construction of the directional hierarchy with (ii) an entry-efficient and robust construction of the nested representation based on the \ac{IACA}. 
For (i), we choose the number of directions for each cluster based on its electrical size, which decouples the number of directions from the level of the tree.
Because recursive face refinement generally cannot produce the resulting numbers of directions, we use spherical-Fibonacci point sets, which remain nearly uniform on the sphere for arbitrary set sizes \cite{keinert_SphericalFibonacciMapping_2015}.
An angular nearest-neighbor mapping links the sets of directions into a directional hierarchy and defines the parent--child relations required by the nested directional representation.
For (ii), we observe that block structures can arise that jeopardize the robustness of the \ac{IACA} and thus an accurate compression, similar to the \ac{MFIE} \cite{tetzner_AdaptiveCrossApproximation_2024}.
The block structures arise in the vector-potential part of the \ac{EFIE}, whose entries depend on the spatial orientation of the basis functions.
We robustify the tree-mimicry pivoting of the \ac{IACA} by enforcing that basis functions of different orientations must be sampled.
The numerical results demonstrate the efficiency and robustness of our approach both for canonical and realistic geometries, exhibiting the expected $\mc O(N\log N)$ scaling for high-frequency problems (i.e., where the frequency is increased in lockstep with $N$) while maintaining the prescribed accuracy.
Preliminary results have been presented at conferences \cite{tetzner_Directional$mathcalH^2$matrix_2026,tetzner_OrientationAwarePivotingr_2026}.

The remainder of this article is organized as follows.
\Cref{sec:Notation} fixes the notation for the \ac{EFIE}, the cluster trees, and the low- and high-frequency admissibility conditions.
\Cref{sec:DirInteractionHierarchy} constructs the geometry-adaptive directions and resulting hierarchy of directionally subdivided interactions associated with one cluster, and \Cref{sec:bottomupWNCA} develops the directional representation using the orientation-aware \ac{IACA}.
\Cref{sec:NumericalResults} reports the numerical validation.

\section{Background and Notation}\label{sec:Notation}
When a time-harmonic wave $(\veg e^\mr{i}, \veg h^\mr{i})$ with electric field $\veg e^\mr{i}$ and magnetic field $\veg h^\mr{i}$ impinges on a perfectly electrically conducting surface $\Gamma$, it induces a surface current density $\veg j$, from which the radiation potentials recover the scattered wave~$(\veg e^\mr{s}, \veg h^\mr{s})$.
We obtain $\veg j$ by solving the \ac{EFIE} for the exterior problem \cite{maue_ZurFormulierungAllgemeinen_1949},
\begin{equation}
  \vecop T \veg j = \(\jm \kappa \TA + \jm \kappa^{-1}\TPhi\)\veg j = \n \times \veg e^\mr{i}\,,
  \label{eq:efie_continuous}
\end{equation}
where $\kappa$ is the wavenumber,
\begin{equation}
  (\TA \veg j)(\veg r) = \n(\veg r) \times \int_\Gamma G(\veg r, \veg r')\veg j(\veg r')\,\dd S(\veg r')
\end{equation}
is the vector-potential contribution, and
\begin{equation}
  (\TPhi \veg j)(\veg r) = \n(\veg r) \times \nabla_\Gamma \int_\Gamma G(\veg r, \veg r')\nabla_\Gamma' \cdot \veg j(\veg r')\,\dd S(\veg r')
\end{equation}
is the scalar-potential contribution.
Here, $\n(\veg r)$ denotes the surface normal at the observation point, and $\nabla_\Gamma$ and $\nabla_\Gamma'$ the surface differential operators with respect to $\veg r$ and $\veg r'$; the free-space Green's function is
\begin{equation}
  G(\veg r, \veg r') = \frac{\e^{-\jm \kappa \abs{\veg r - \veg r'}}}{4 \uppi \abs{\veg r - \veg r'}}\,.
\end{equation}
We discretize \eqref{eq:efie_continuous} using the \ac{MoM}~\cite{harrington_FieldComputationMoment_1993} with a Petrov--Galerkin ansatz, using \ac{RWG} basis functions $\veg \fk_\beta$ for expansion and $\n \times \veg \fk_\alpha$ for testing \cite{rao_ElectromagneticScatteringSurfaces_1982}, which yields the linear system
\begin{equation}
  \mat T \vec j \equiv \(\jm \kappa \matTA + \jm \kappa^{-1}\matTPhi\)\vec j = \vec e^\mr{i}\,,
  \label{eq:efie_discrete}
\end{equation}
with
\begin{equation}
  \begin{aligned}
    \sbr{\matTA}_{\alpha\beta}   & = \duality{\n \times \veg \fk_\alpha}{\TA \veg \fk_\beta}\,,      \\
    \sbr{\matTPhi}_{\alpha\beta} & = \duality{\n \times \veg \fk_\alpha}{\TPhi \veg \fk_\beta}\,,    \\
    \sbr{\vec e^\mr{i}}_{\alpha} & = \duality{\n \times \veg \fk_\alpha}{\n \times \veg e^\mr{i}}\,,
  \end{aligned}
  \label{eq:efie_matrix_entries}
\end{equation}
where the $L^2$-inner product on $\Gamma$ is defined as
\begin{equation}
  \duality{\veg f}{\veg g}
  \coloneqq
  \int_\Gamma
  \veg f(\veg r)\cdot\veg g(\veg r)\,\dd S(\veg r)\,.
  \label{eq:surface_pairing}
\end{equation}

The matrix $\mat T$ is dense, but suitable sub-blocks admit low-rank approximations.
To expose them, we collect the indices of the testing and expansion functions in the sets $I$ and $J$, respectively, and organize these sets into cluster trees $T_I$ and $T_J$~\cite[Sec.~5.3]{hackbusch_HierarchicalMatricesAlgorithms_2015}.
Each cluster $t \in T_I$ is a subset of $I$, with the root holding all row indices.
Starting from the root, we recursively split each cluster into children until it has fewer than $n_\mr{min}$ indices, and we denote the children of $t$ by $S(t)$.
A cluster $t$ is a leaf if $S(t)=\emptyset$; otherwise, it is the disjoint union of its children,
\begin{equation}
  t = \dot\bigcup_{t' \in S(t)} t' \, .
  \label{eq:cluster_tree_partition}
\end{equation}
We build the column tree $T_J$ analogously.

A sub-block $\mat T_{t,s}$ describes the interaction between a row cluster $t \in T_I$ and a column cluster $s \in T_J$.
We denote the corresponding geometrical supports, given by the unions of the supports of the associated basis functions, by $X_t$ and $Y_s$, respectively.
For bounded sets $X$ and $Y$, we define
\begin{align}
  \diam(X)   & \coloneqq \sup_{\veg x,\veg y\in X}\abs{\veg x-\veg y}\,, \\
  \dist(X,Y) & \coloneqq \inf_{\substack{\veg x\in X                   \\ \veg y\in Y}}\abs{\veg x - \veg y}\,.
\end{align}
The mathematical conditions are stated in terms of the diameters of the supports and the distance between them.
In the implementation, these quantities are estimated using bounding balls $B_t$ and $B_s$ with centers $\veg c_t$ and $\veg c_s$ and radii $ r_t$ and $r_s$, respectively:
For a K-means tree \cite[Sec. 9.3]{gibson_MethodMomentsElectromagnetics_2021}, \cite{fukunaga_BranchBoundAlgorithm_1975}, $\veg c_t$ is the K-means center and $B_t$ is the smallest ball centered at $\veg c_t$ that encloses $X_t$ so that $\diam (B_t) \geq\diam (X_t)$.
For an octree, $\veg c_t$ is the box center and $B_t$ is the ball enclosing the box, whose diameter is $\sqrt{3}$ times the side length of the box.
Accordingly, the implementation estimates $\diam(X_t)$ by $\diam(B_t)=2r_t$ and $\dist(X_t,Y_s)$ by
\begin{equation}
  \dist(B_t,B_s)
  =
  \max\{\abs{\veg c_t-\veg c_s}-r_t-r_s,0\}\,.
\end{equation}
The quantities for column clusters are defined analogously.

We obtain the admissible blocks by recursively subdividing the root block $(I,J)$ until each block is admissible or cannot be subdivided further.
The admissibility criterion depends on the electrical size of each of the interacting clusters.
Following \cite{bebendorf_WidebandNestedCross_2015,borm_Directional$mathcalH^2$matrix_2017}, if $\kappa \min\{\diam (X_t),\diam (Y_s)\} \leq 1$, the block is a low-frequency interaction and is admissible when
\begin{equation}
  \eta_{\mr l}\, \dist(X_t, Y_s) \geq \max\{\diam (X_t), \diam (Y_s)\}\,.
  \label{eq:lfadm}
\end{equation}
If $\kappa \min\{\diam (X_t),\diam (Y_s)\} > 1$, it is a high-frequency interaction and is admissible when
\begin{equation}
  \eta_{\mr h}\, \dist(X_t, Y_s) \geq \kappa \max\{\diam (X_t)^2,\diam( Y_s)^2\}\,.
  \label{eq:hfadm}
\end{equation}
Here $\eta_{\mr l},\eta_{\mr h} > 0$ are admissibility parameters.

This subdivision yields a partition $P$ of $I \times J$ into admissible and non-admissible blocks, that is,
\begin{equation}
  P = P_\mr{adm}\, \dot\cup \, P_\mr{nonadm}\,.
\end{equation}
The admissible blocks are further separated according to the frequency regime, that is,
\begin{equation}
  P_\mr{adm} = P_\mr{low} \,\dot\cup \, P_\mr{high}\,,
\end{equation}
where $P_\mr{low}$ and $P_\mr{high}$ contain the admissible low- and high-frequency interactions, respectively.
The blocks in $P_\mr{nonadm}$ are stored explicitly (i.e., uncompressed).

Since an admissible interaction is not subdivided further, each pair $(t,s)\in P_\mr{adm}$ is encountered at the first stage of the recursive subdivision at which the clusters $t$ and $s$ are referred to as well separated.
For each row cluster $t \in T_I$, we collect the indices of the well-separated column clusters in
\begin{equation}
  \mc F_\mr{w}(t)
  \coloneqq
  \bigcup
  \{s \in T_J : (t,s) \in P_\mr{adm}\}
  \subseteq J\,.
  \label{eq:well_separated_clusters}
\end{equation}

The low-frequency interactions in $P_\mr{low}$ admit low-rank approximations with ranks bounded independently of the wavenumber, and we represent them by standard nested low-rank approximations constructed with the \ac{NCA} \cite{bebendorf_ConstructingNestedBases_2012}.
For the high-frequency interactions in $P_\mr{high}$, the oscillatory factor $\e^{-\jm\kappa\abs{\veg r-\veg r'}}$ of the free-space Green's function causes the rank required for a prescribed accuracy to grow with the electrical size of the interacting supports \cite{brick_InterpretingMomentMatrix_2026, bebendorf_WidebandNestedCross_2015}.

The low- or high-frequency classification, \eqref{eq:lfadm} and \eqref{eq:hfadm}, determines which admissibility condition is applied to an interaction $(t, s)$.
A cluster $t$ generally participates in several admissible interactions, but its basis is constructed once and reused for all of them.
It must therefore carry either a single non-directional basis in the low-frequency regime or a family of directional bases in the high-frequency regime, independently of the well-separated clusters $s$.
The same cluster regime determines the transfer matrices that connect the basis of a parent cluster to those of its children.
Since one such nested relation involves all children of a parent, we assign the regime jointly to all siblings.
Specifically, we use the average of their bounding-ball diameters and classify all siblings as low-frequency clusters if $\kappa$ times this average does not exceed one; otherwise, we classify them as high-frequency clusters.

For a low-frequency cluster $t$, the set $\mc F(t)$ contains $\mc F_\mr{w}(t)$ together with the indices inherited from its parent $t_\mr{p}$ if $t_\mr{p}$ is also low frequency.
Accordingly, we define
\begin{equation}
  \mc F(t)
  \coloneqq
  \begin{cases}
    \mc F_\mr{w}(t) \cup \mc F(t_\mr p)\,,
    & \text{if $t_\mr p$ is low frequency}, \\
    \mc F_\mr{w}(t)\,,
    & \text{otherwise.}
  \end{cases}
  \label{eq:well_separated_cluster_inheritance}
\end{equation}

For a high-frequency cluster $t$, the oscillatory factor of Green's function $\e^{-\jm\kappa\abs{\veg r-\veg r'}}$ can be approximated by the plane wave $\e^{-\jm\kappa\langle\veg r-\veg r',\veg e\rangle}$ if the point $\veg r'$ is located in a well-separated cluster $s$ within a sufficiently narrow cone around a propagation direction $\veg e$, where $\langle\cdot,\cdot\rangle$ denotes the Euclidean inner product.
Factoring out this plane wave leaves a slowly varying kernel with bounded rank, which can be compressed by nested low-rank approximations \cite{bebendorf_WidebandNestedCross_2015,borm_Directional$mathcalH^2$matrix_2017}.
Following \cite{bebendorf_WidebandNestedCross_2015}, the number of required directions is dictated by the condition that, for every $s\subseteq\mc F_\mr{w}(t)$, at least one $\veg e\in\mc E(t)$ satisfies
\begin{equation}
  \abs{\sin \measuredangle(\veg c_t-\veg y,\veg e)}
  \leq
  \frac{\gamma}{\kappa \diam (X_t)}\,,
  \qquad
  \veg y \in Y_s \,,
  \label{eq:anglecondition}
\end{equation}
where $\gamma>0$ controls the directional resolution \cite{bebendorf_WidebandNestedCross_2015}.
With this definition, the vectors from the source points $\veg y\in Y_s$ to the row-cluster center $\veg c_t$ lie within a thin cone around $\veg e$, and the cone narrows as the electrical size $\kappa\diam (X_t)$ grows, so electrically larger clusters need more directions to be covered.

For a high-frequency cluster $t$, we associate each well-separated cluster $s$ with a direction $\veg e (t, s)\in \mc E(t)$.
These directions partition the well-separated clusters into
\begin{equation}
  \mc F_{\veg e,\mr{w}}(t)
  \coloneqq
  \bigcup
  \{s \in T_J : (t,s) \in P_\mr{adm}\,,\ \veg e(t,s)=\veg e\} \subseteq J\,.
  \label{eq:directional_well_separated_clusters}
\end{equation}
The corresponding quantities for column clusters are defined analogously with rows and columns interchanged; in particular, $\veg e(s,t)\in\mc E(s)$ denotes the direction associated with the interaction $(t,s)$ from the column-cluster side.

\begin{figure}[t]
  \centering
  \includegraphics[width=0.93\columnwidth]{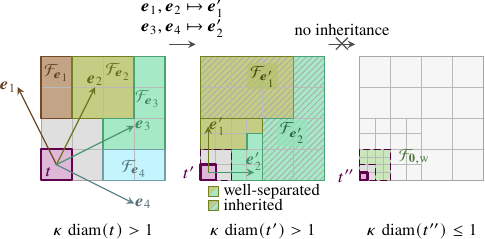}
  \caption{Inheritance of directionally subdivided and well separated clusters.
    Solid regions denote well-separated clusters assigned directly to the respective cluster, and hatching denotes inherited clusters; inheritance stops at the high-to-low-frequency transition.}
  \label{fig:directional_interaction_inheritance}
\end{figure}

The sets $\mc F_{\veg e,\mr{w}}(t)$ contain the indices of the well-separated clusters assigned directly to $t$.
We denote by $\mc F_{\veg e}(t)$ the corresponding directionally subdivided sets including those inherited from parent clusters, as illustrated in \Cref{fig:directional_interaction_inheritance}.
As mentioned, inheritance is restricted to parent--child links within the same frequency regime.
Accordingly, a low-frequency cluster inherits only from a low-frequency parent, while a low-frequency cluster with a high-frequency parent starts from $\mc F_\mr{w}(t)$.
A mathematical definition of $\mc F_{\veg e}(t)$ follows in \Cref{sec:DirInteractionHierarchy} after introducing the parent--child mapping of the directions.

For an admissible interaction $(t,s)$ of two high-frequency clusters $t$ and $s$ with $s\subseteq\mc F_{\veg e,\mr{w}}(t)$, the associated sub-block admits the nested low-rank approximation
\begin{equation}
  \mat T_{t,s}
  \approx
  \mat U_t^{\veg e(t, s)}\,\mat S_{t,s}\,\mat V_s^{\veg e(s, t)}\,,
  \label{eq:dh2_factorization}
\end{equation}
where $\mat S_{t,s}$ is a coupling matrix and $\mat U_t^{\veg e(t, s)}$, $\mat V_s^{\veg e(s, t)}$ are the directional cluster bases associated with the respective direction.
For a low-frequency cluster, the corresponding directional superscript is omitted and the non-directional basis $\mat U_t$ or $\mat V_s$ is used instead.
For mixed-regime interactions, the directional superscript is omitted for the low-frequency cluster.

Within each regime, transfer matrices represent the bases in nested form and coupling matrices are stored per admissible interaction \cite{bebendorf_ConstructingNestedBases_2012,bebendorf_WidebandNestedCross_2015}; no transfer connects the two regimes.
To realize this representation on geometry-adaptive cluster trees such as a K-means tree, the sets of directions must be chosen according to the local cluster geometry and connected explicitly along the cluster tree.

\section{Geometry-Adaptive Direction Sets and Directional Hierarchy}
\label{sec:DirInteractionHierarchy}

\begin{figure}[]
  \centering
  \includegraphics[width=0.85\columnwidth]{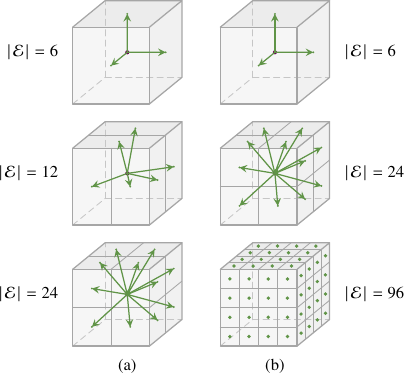}
  \caption{Cube base directional subdivision:
    (a) successive binary refinement \cite{bebendorf_WidebandNestedCross_2015};
    (b) successive four-way refinement \cite{borm_Directional$mathcalH^2$matrix_2017}.
    Labels give the total number of directions over all six faces; for $\lvert\mc E\rvert=96$, only direction endpoints are shown.}
  \label{fig:prescribed_direction_refinements}
  \phantomsubcaption\label{fig:prescribed_direction_refinements_bebendorf}
  \phantomsubcaption\label{fig:prescribed_direction_refinements_boerm}
\end{figure}

Existing \ac{DH2} constructions \cite{bebendorf_WidebandNestedCross_2015,borm_Directional$mathcalH^2$matrix_2017} obtain the sets of directions and their parent--child relations from prescribed angular refinements.
In \cite{bebendorf_WidebandNestedCross_2015}, each cluster $t$ is associated with a hierarchy of pyramidal cones that defines its directional subdivision and the corresponding parent--child mapping.
The construction starts from the centers of the six faces of a cube and recursively splits each face region perpendicular to its longest edge.
The centers of the resulting face regions define the directions, while the refinement hierarchy defines their parent--child mapping, as illustrated in \Cref{fig:prescribed_direction_refinements_bebendorf}.

In \cite{borm_Directional$mathcalH^2$matrix_2017}, the subdivision is also based on cube faces but follows the levels of an octree.
Starting from the six face centers, each refinement subdivides every face region into four squares and projects their centers onto the unit sphere.
Each refinement therefore quadruples the number of directions, matching the quadratic growth of the required directional resolution when the cluster diameter doubles, as shown in \Cref{fig:prescribed_direction_refinements_boerm}.
The $\lvert\mc E\rvert=24$ stages coincide: two binary splits of \cite{bebendorf_WidebandNestedCross_2015} produce the same four-way subdivision as one four-square split of \cite{borm_Directional$mathcalH^2$matrix_2017}.

Both constructions associate each tree level with a fixed set of directions.
This is natural for regularly refined box trees, where the cluster diameter is related to the level, but becomes restrictive for geometry-adaptive trees with clusters of different diameters on the same level.
Moreover, recursive face refinement restricts the number of directions to values generated by repeated doubling~\cite{bebendorf_WidebandNestedCross_2015} or quadrupling~\cite{borm_Directional$mathcalH^2$matrix_2017}.
In contrast, the upper bound for the angular resolution \eqref{eq:anglecondition} depends continuously on the individual cluster diameter and therefore generally leads to a cluster-specific number of directions that does not coincide with these prescribed refinement stages.
A common level-wise set of directions then over-resolves smaller clusters and increases storage and \ac{MVP} cost.
We therefore determine the number of directions for each cluster from its electrical size $\kappa\diam(B_t)$ and construct the corresponding set of directions independently.
Since this construction requires only the bounding ball $B_t$ defined in \Cref{sec:Notation}, it is applicable independently of the underlying cluster-tree construction, provided that for each cluster $t$ we have a bounding ball $B_t$.

\subsection{Spherical-Fibonacci Direction Sets}

To derive the adaptive sets of directions, we first consider the planar case, where both the placement of the directions and their angular covering can be characterized exactly.
For a planar geometry, thus essentially two-dimensional, the directions are restricted to the geometry plane and
are therefore vectors from the coordinate origin covering the unit circle.
For example, a construction for planar antenna arrays is presented in~\cite{sewraj_MacroBasisFunction_2024}, where the azimuthal directional subdivision is adapted to the electrical size of the associated cluster.
In the planar setting, equally spaced directions provide a uniform covering of the unit circle.
For $n_t$ directions, we define
\begin{equation}
  \veg e_k = \(\cos(2\uppi k/n_t),\, \sin(2\uppi k/n_t)\)\,,\qquad k = 0,\ldots,n_t-1 \,,
  \label{eq:dir_2d_points}
\end{equation}
and set $\mc E(t)=\{\veg e_k : k=0,\ldots,n_t-1\}$.
Every unit vector then has an angular distance of at most $\uppi/n_t$ to its closest direction.
Consequently, matching the angular resolution prescribed by \eqref{eq:anglecondition} gives
\begin{equation}
  \sin\(\frac{\uppi}{n_t}\) \leq \frac{\gamma}{\kappa\,\diam(B_t)}\,,
  \label{eq:dir_2d_condition}
\end{equation}
which yields
\begin{equation}
  n_t \geq \frac{\uppi}{\arcsin(\gamma/(\kappa\,\diam(B_t)))} \,.
  \label{eq:dir_2d_count}
\end{equation}
Hence, the electrical size of $t$ determines the number of directions independently of the level of $t$ in the cluster tree.

Analogously, for three-dimensional geometries, the directions are vectors from the coordinate origin covering the unit sphere.
Unlike the planar case, however, no uniformly spaced point sets with an exactly known covering angle exist for arbitrary numbers of directions.
We therefore first derive an estimate for the number of directions required for $\kappa \diam(B_t)$ and subsequently place this number of directions approximately uniformly on the sphere.
To obtain such an estimate, we start by considering the six initial directions employed in the established directional subdivisions~\cite{bebendorf_WidebandNestedCross_2015,borm_Directional$mathcalH^2$matrix_2017} (see $\abs{\mc E} = 6$ in \Cref{fig:prescribed_direction_refinements}), which are a uniform covering of the sphere.
These directions point toward the centers of the six faces of a cube and correspond to the six coordinate directions, or equivalently to the vertices of an octahedron.
Their maximal angular distance to any point on the unit sphere is
\begin{equation}
  \alpha_0 = \arccos \(\frac{1}{\sqrt{3}}\)\,.
  \label{eq:octahedron_initial_angle}
\end{equation}
For $n_t$ approximately uniformly distributed directions, the spherical area assigned to each direction scales as $1/n_t$.
We assume for the sake of simplicity that each area covering the sphere corresponds to the base area of a cone from the origin.
Then the maximal angular distance $\alpha$ corresponds to the half-angle of the cones.
Clearly, $\alpha$ is proportional to the square root of their base area and thus
\begin{equation}
  \alpha(n_t)\approx\alpha_0\sqrt{\frac{6}{n_t}} \,.
  \label{eq:angular_scaling}
\end{equation}
To satisfy the angular resolution dictated by \eqref{eq:anglecondition}, we require the estimated covering angle to meet
\begin{equation}
  \sin\(\alpha(n_t)\)
  \leq
  \frac{\gamma}{\kappa\diam(B_t)}\,.
\end{equation}
Inserting \eqref{eq:angular_scaling}, we find
\begin{equation}
  n_t =\left\lceil 6 \(\frac{\alpha_0}{\arcsin(\gamma/(\kappa\,\diam(B_t)))}\)^2\right\rceil \,.
  \label{eq:dir_3d_count}
\end{equation}

Because the regime is assigned from the average bounding-ball diameter of two siblings, a cluster $t$ might be low-frequency, that is, $\kappa\diam(B_t)\leq1$.
In this case, if $\gamma=1$ the argument of the $\arcsin$ exceeds one and we thus limit it to one, yielding the coarsest possible set of directions.

Equation \eqref{eq:dir_3d_count} determines only the necessary number of directions.
To realize the corresponding directions providing a (near-)uniform covering of the unit sphere, we use spherical-Fibonacci points, which are available for arbitrary values of $n_t$~\cite{keinert_SphericalFibonacciMapping_2015}.
For $i=0,\ldots,n_t-1$, let
\begin{align}
  z_i &= 1 - \frac{2(i+0.5)}{n_t}\, ,\\  r_i &= \sqrt{1-z_i^2} \, , \\ \phi &= \uppi(3-\sqrt{5}) \, ,
  \label{eq:fibonacci_parameters}
\end{align}
and the Fibonacci points are defined by
\begin{equation}
  \veg e_i = \(r_i\cos(i\phi),\, r_i\sin(i\phi),\, z_i\)\,.
  \label{eq:fibonacci_direction}
\end{equation}
Then $\mc E(t)=\{\veg e_i : i=0,\ldots,n_t-1\}$ is the set of directions for a high-frequency cluster.

In contrast to the unit-circle case, \eqref{eq:dir_3d_count} is an estimate and not a sharp lower bound.
Clearly, over-resolution introduces more directions than necessary, while under-resolution assigns a wider angular range to each direction and generally increases the approximation ranks.
For severe under-resolution, the bases of electrically smaller child clusters may no longer represent the higher-rank interactions required at their parents with sufficient accuracy, so that additional errors are introduced through the nested representation.
Such an effect would become visible as an increasing deviation of the approximation error from the desired compression accuracy for electrically larger problems.
In our numerical experiments, we do not observe such a failure, indicating the applicability of our heuristic.
In fact, $\diam(B_t)\geq\diam (X_t)$ is a conservative estimate, biasing the number of directions toward finer resolution.

The resulting set of directions $\mc E(t)$ defines the directional subdivision of the admissible interactions of $t$.
For a high-frequency cluster $t$ and $s \subseteq \mc F_\mr{w}(t)$, we assign $s$ to its closest direction, that is,
\begin{equation}
  \veg e(t,s) = \arg\min_{\veg e \in \mc E(t)} \measuredangle(\veg c_t-\veg c_s,\veg e)\, ,
  \label{eq:direction_assignment}
\end{equation}
which defines the directional subsets $\mc F_{\veg e,\mr{w}}(t)$ from \eqref{eq:directional_well_separated_clusters}.

We next link the directions of $t$ to the directions of its children $S(t)$.
Unlike the subdivisions in \cite{bebendorf_WidebandNestedCross_2015, borm_Directional$mathcalH^2$matrix_2017}, each cluster $t$ has an individual set of directions, and we therefore need to map each direction of $t$ to one direction of every child $t' \in S(t)$.

\subsection{Directional Hierarchy}

\begin{figure}
  \centering
  \includegraphics{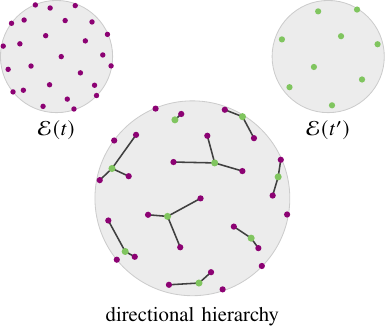}
  \caption{Directional hierarchy for locally generated sets of directions. Each parent direction is mapped to the direction of the child that is closest in angle, and several parent directions may be mapped to the same child direction.}
  \label{fig:directional_hierarchy}
\end{figure}

We use a nearest-neighbor mapping with respect to the angle between two directions.
Specifically, if $t' \in S(t)$ be a high-frequency child of a high-frequency parent $t$ with a direction $\veg e \in \mc E(t)$, we map $\veg e$ to the direction $\veg e'$ of $t'$ that is closest in angle, that is,
\begin{equation}
  \veg e'(\veg e,t')
  =
  \arg\min_{\widetilde{\veg e} \in \mc E(t')}
  \measuredangle(\veg e,\widetilde{\veg e}) \,.
  \label{eq:dirmapping}
\end{equation}
In the case of two or more equally distant directions, an arbitrary choice is made.
Since the child cluster is electrically no larger than the parent, it has fewer directions and several parent directions may be mapped to the same child direction.
\Cref{fig:directional_hierarchy} illustrates this many-to-one relation.
For a high-frequency child $t' \in S(t)$ of a high-frequency parent $t$ and a direction $\veg e' \in \mc E(t')$, we define
\begin{equation}
  \mc F_{\veg e'}(t')
  \coloneqq
  \mc F_{\veg e',\mr{w}}(t')
  \cup
  \bigcup_{\substack{\veg e \in \mc E(t)\\
      \veg e'(\veg e,t')=\veg e'}}
  \mc F_{\veg e}(t)\,.
  \label{eq:directional_interaction_inheritance}
\end{equation}
For the root cluster, we set
$\mc F_{\veg e}(t)=\mc F_{\veg e,\mr{w}}(t)$.
We note that this process terminates when the children $t' \in S(t)$ are low-frequency.

\section{Bottom-Up Wideband NCA Construction}
\label{sec:bottomupWNCA}

The directional hierarchy constructed in \Cref{sec:DirInteractionHierarchy} provides $\mc F_{\veg e}(t)$ for every cluster $t$ and direction $\veg e\in\mc E(t)$.
Whenever $\mc F_{\veg e}(t)\neq\emptyset$, we select pivots from $\mat T_{t,\mc F_{\veg e}(t)}$ and construct the corresponding bases and transfer matrices in a bottom-up traversal of the cluster tree, proceeding from the leaves toward the root~\cite{tetzner_IncompleteAdaptiveCross_2026}.
We formulate the construction for the row tree $T_I$; the column-tree construction follows with rows and columns interchanged.

\subsection{Construction of Directional Cluster Bases}
\label{subsec:iaca_wideband_nca}
For low-frequency clusters, the pivots, bases, and transfer matrices are constructed as introduced in \cite{tetzner_IncompleteAdaptiveCross_2026}.
In the following, we consider the corresponding \ac{DH2} representation for high-frequency clusters~\cite{bebendorf_WidebandNestedCross_2015}; for the corresponding low-frequency expressions, the directional superscripts and subscripts are omitted, for example, $\mc F_{\veg e}(t)$ is replaced by $\mc F(t)$.
For each $t \in T_I$ and $\veg e \in \mc E(t)$ with $\mc F_{\veg e}(t) \neq \emptyset$, we select row indices $\tau_t^{\veg e}\subseteq t$ and column indices $\sigma_{\mc F_{\veg e}(t)}\subseteq \mc F_{\veg e}(t)$.
Here, $\sigma_{\mc F_{\veg e}(t)}$ denotes  the set of column indices associated with the basis functions of the clusters $s \subseteq \mc F_{\veg e}(t)$.
These pivot sets define the CUR approximation
\begin{equation}
  \mat T_{t,\mc F_{\veg e}(t)}
  \approx
  \mat T_{t,\sigma_{\mc F_{\veg e}(t)}}
  \(\mat T_{\tau_t^{\veg e},\sigma_{\mc F_{\veg e}(t)}}\)^{-1}
  \mat T_{\tau_t^{\veg e},\mc F_{\veg e}(t)} \,.
  \label{eq:directional_cur}
\end{equation}

We determine these pivot sets with the \ac{IACA}~\cite{tetzner_IncompleteAdaptiveCross_2026}.
Compared with a standard \ac{ACA}, the \ac{IACA} selects the column pivots $\sigma_{\mc F_{\veg e}(t)}$ geometrically and thereby avoids assembling complete rows of $\mat T_{t,\mc F_{\veg e}(t)}$, resulting in a faster construction of the \ac{DH2} representation.
For each selected column, subtracting its approximation obtained from the previously selected pivots yields a residual column, whose maximum absolute value determines the row pivot in $\tau_t^{\veg e}$.
Since the \ac{IACA} does not assemble the complete matrix $\mat T_{\tau_t^{\veg e},\mc F_{\veg e}(t)}$ during pivoting, only the sampled residual columns are used to estimate the Frobenius norm of the residual matrix
\begin{equation}
  \normF{\wmat T_{t,\mc F_{\veg e}(t)}} = \normF{
    \mat T_{t,\mc F_{\veg e}(t)}
    -
    \mat T_{t,\sigma_{\mc F_{\veg e}(t)}}
    \(\mat T_{\tau_t^{\veg e},\sigma_{\mc F_{\veg e}(t)}}\)^{-1}
    \mat T_{\tau_t^{\veg e},\mc F_{\veg e}(t)} }\,.
  \label{eq:residual_matrix}
\end{equation}

With the computed pivots, we form the row basis
\begin{equation}
  \mat U_t^{\veg e}
  =
  \mat T_{t,\sigma_{\mc F_{\veg e}(t)}}
  \(\mat T_{\tau_t^{\veg e},\sigma_{\mc F_{\veg e}(t)}}\)^{-1}
  \label{eq:directional_basis}
\end{equation}
if $t$ is a leaf or if $t$ is a high-frequency cluster whose children are low-frequency.
Otherwise, $t$ and its children share the same frequency regime and we construct transfer matrices.
Given a direction $\veg e \in \mc E(t)$, for each child $t'_i\in S(t)$, let $\veg e'_i\coloneqq\veg e'(\veg e,t'_i)$ be the corresponding child direction.
The transfer matrix is
\begin{equation}
  \mat \Theta^{\veg e}_{t'_i,t}
  =
  \mat T_{\tau_{t'_i}^{\veg e'_i},\sigma_{\mc F_{\veg e}(t)}}
  \(\mat T_{\tau_t^{\veg e},\sigma_{\mc F_{\veg e}(t)}}\)^{-1}\,.
  \label{eq:directional_transfer}
\end{equation}
The first matrix factor is a submatrix of $\mat T_{\tilde t, \mc F_{\veg e}(t)}$ with $\tilde t$ defined by \cite[Eq. (18)]{tetzner_IncompleteAdaptiveCross_2026}.
The resulting transfer matrices represent the parent basis in nested form, that is, 
\begin{equation}
  \mat U_t^{\veg e}
  =
  \begin{pmatrix}
    \mat U_{t'_1}^{\veg e'_1} & & 0 \\
    & \ddots & \\
    0 & & \mat U_{t'_{\abs{S(t)}}}^{\veg e'_{\abs{S(t)}}}
  \end{pmatrix}
  \begin{pmatrix}
    \mat \Theta_{t'_1,t}^{\veg e} \\
    \vdots \\
    \mat \Theta_{t'_{\abs{S(t)}},t}^{\veg e}
  \end{pmatrix}\,.
  \label{eq:directional_nestedness}
\end{equation}

Analogously, the column-tree traversal constructs $\mat V^{\veg e}_s$ and $\sigma_s^{\veg e}$ for the directions $\veg e\in \mc E(s)$.
An admissible block $(t,s)$ of two high-frequency clusters is then represented by \eqref{eq:dh2_factorization} with the coupling matrix
\begin{equation}
  \mat S_{t,s}
  =
  \mat T_{\tau_t^{\veg e(t, s)},\sigma_s^{\veg e(s, t)}}\,.
  \label{eq:directional_coupling}
\end{equation}


\subsection{Orientation-Aware Pivoting}
\label{subsec:orientation_aware_pivoting}

The accuracy of the \ac{IACA} construction depends on the representative selection of the column pivots $\sigma_{\mc F_{\veg e}(t)}$.
It is well known that the \ac{ACA} may suffer from premature convergence~\cite{heldring_ConvergenceACAAlgorithm_2014,heldring_ImprovingAccuracyAdaptive_2021,zhang_ImprovedAdaptiveCross_2025,almuna-morales_RobustHierarchicalMatrix_2026a} and, clearly, methods derived from the \ac{ACA}, such as the \ac{IACA}, can inherit this behavior.
In the case of the \ac{MFIE}, we identified block-structured interaction matrices as a possible source of premature convergence \cite{tetzner_AdaptiveCrossApproximation_2024}.
Such block structures are, for example, likely to occur whenever the underlying mesh is structured so that only a few orientations of expansion and testing functions are possible.
Related difficulties in obtaining representative low-rank approximations of \ac{MoM} matrix blocks with vector basis and testing functions have been reported in~\cite{kelley_RobustGridBasedComputation_2023}.

\begin{figure}
  \centering
  \begin{minipage}{0.49\linewidth}
    \centering
    \includegraphics{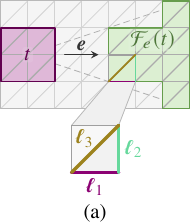}
  \end{minipage}
  \hfill
  \begin{minipage}{0.49\linewidth}
    \centering
    \includegraphics{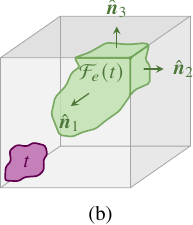}
  \end{minipage}
  \caption{Two geometrical configurations leading to orientation-dependent coupling within $\mc F_{\veg e}(t)$. (a) $\mc F_{\veg e}(t)$ associated with direction $\veg e$ contains \acp{RWG} with the different edge vectors $\veg\ell_1$, $\veg\ell_2$, and $\veg\ell_3$. (b) $\mc F_{\veg e}(t)$ extends over the three faces adjacent to a corner of the cube, forming three nearly planar groups with the orthogonal surface normals $\hat{\veg n}_1$, $\hat{\veg n}_2$, and $\hat{\veg n}_3$.}
  \label{fig:orientation_failures}
  \phantomsubcaption\label{fig:orientation_failures_edges}
  \phantomsubcaption\label{fig:orientation_failures_normals}
\end{figure}

Similar to the \ac{MFIE} operator, the discretization of $\TA$ depends on the orientation of the expansion and testing functions $\veg\fk_\beta$ and $\n \times \veg \fk_\alpha$.
This is in contrast to $\TPhi$, where, in its weak form, the divergence operators eliminate the effect of the orientation.
In fact, for low-frequency problems, the $\TPhi$ contribution mitigates the effect of the orientation-dependency of $\TA$ as $\matTPhi$ dominates $\mat T$ due to the inverse frequency scaling in \eqref{eq:efie_continuous}.
Consequently, no premature convergence was observed in \cite{tetzner_IncompleteAdaptiveCross_2026}, because only low-frequency problems were studied.
Specifically, we find for $\TA$ two configurations that lead to premature convergence, depicted in \Cref{fig:orientation_failures}.

The first type arises for planar domains with a structured mesh such as shown in \Cref{fig:orientation_failures_edges}, affecting both \ac{ACA} and \ac{IACA}.
In this example, the \acp{RWG} are associated with only three edge vectors $\veg\ell_1$, $\veg\ell_2$, and $\veg\ell_3$.
If the expansion and testing functions are sorted according to $\veg \ell_i$, the interaction matrix $\mat T_{\mr A, t,\mc F_{\veg e}(t)}$ becomes a $3 \times 3$ block structure of the form
\begin{equation}
  \mat T_{\mr A, t, \mc F_{\veg e}(t)} = \begin{bmatrix}
    \mat T_{\veg \ell_1, \veg \ell_1} & \mat T_{\veg \ell_1, \veg \ell_2} & \mat T_{\veg \ell_1, \veg \ell_3} \\
    \mat T_{\veg \ell_2, \veg \ell_1} & \mat T_{\veg \ell_2, \veg \ell_2} & \mat T_{\veg \ell_2, \veg \ell_3} \\
    \mat T_{\veg \ell_3, \veg \ell_1} & \mat T_{\veg \ell_3, \veg \ell_2} & \mat T_{\veg \ell_3, \veg \ell_3} \\
  \end{bmatrix}\,.
\end{equation}
Numerically, we observe that the diagonal blocks $\mat T_{\veg\ell_i,\veg\ell_i}$ have larger Frobenius norms than the off-diagonal blocks $\mat T_{\veg\ell_i,\veg\ell_j}$, $i\neq j$, and that this relation persists for the corresponding blocks of the residual matrix during the \ac{ACA} iterations.
If the pivots remain associated with a single $\veg\ell_i$, the corresponding diagonal block and the off-diagonal blocks involving $\veg \ell_i$ are accurately approximated, while the diagonal blocks associated with the remaining $\veg \ell_j$, $j \neq i$, are poorly approximated, leading to larger than desired errors.
Consequently, for commonly used compression tolerances, the \ac{ACA} may converge prematurely before switching to another edge vector $\veg \ell_j$.
A similar problem can occur for the \ac{IACA} if mimicry pivoting, by chance, repeatedly selects expansion functions associated with the same edge vector $\veg \ell_i$.
Then, from the resulting residual columns, analogously to the \ac{ACA}, testing functions associated with the same edge vector $\veg \ell_i$ are chosen since their absolute values will be larger than for the off-diagonal blocks.

A second type arises for a cluster $t$ interacting with two or more orthogonal planar patches contained in  $\mc F_{\veg e}(t)$, for example, the three faces adjacent to a corner of the cube, such as depicted in \Cref{fig:orientation_failures_normals}.
We associate with each patch in $\mc F_{\veg e}(t)$ a surface normal $\hat{\veg n}_i$.
Neglecting basis functions located on the edges connecting the patches, we can group the remaining basis functions according to the patch that they are defined on.
Thus, for the example, we obtain a $1 \times 3$ block structure
\begin{equation}
  \matTA[t,\mc F_{\veg e}(t)]
  \approx
  \begin{bmatrix}
    \mat T_{\hat{\veg n}_1} & \mat T_{\hat{\veg n}_2} & \mat T_{\hat{\veg n}_3}
  \end{bmatrix}\,,
  \label{eq:orientation_groups}
\end{equation}
Numerically, we observe that during the factorization the residual norm of one of these blocks may already be sufficiently small to meet the prescribed tolerance while the other blocks are still poorly approximated.
For the \ac{ACA}, this is generally uncritical because each selected residual row contains entries from all three blocks and can therefore expose approximation errors in any of them.
For the \ac{IACA}, the convergence criterion (see \cite{tetzner_IncompleteAdaptiveCross_2026} for details) is 
\begin{equation}
  \norm{\vec a_r}_F
  \leq
  \varepsilon\,
  \frac{\normF{\mat T_{t,\sigma^{(r)}_{\mc F_{\veg e}(t)}}}}{\sqrt{r}}\,,
  \label{eq:iaca_convergence}
\end{equation}
with $\vec a_r$ denoting the sampled column from $\wmat T_{t, \mc F_{\veg e}(t)}$ at step $r$ and $\sigma^{(r)}_{\mc F_{\veg e}(t)}$ the column indices selected up to this step.
The column $\veg a_r$ belongs to a single block $\mat T_{\hat{\veg n}_i}$ and its norm may therefore already be small although the residuals of the blocks associated with $\hat{\veg n}_j$, $j\neq i$, remain large.
Consequently, the convergence criterion of the \ac{IACA} may be satisfied while significant approximation errors in $\mat T_{\hat{\veg n}_j}$, $j\neq i$ remain undetected.

We address both configurations by extending the tree mimicry pivoting of the \ac{IACA}~\cite{tetzner_IncompleteAdaptiveCross_2026}.
Let $Z_\sigma$ contain the positions of the previously selected \acp{RWG} $\sigma_{\mc F_{\veg e}(t)}$, where we define the position of an \ac{RWG} as the center of its associated edge.
For $r>1$, we define
\begin{equation}
  \rho_r(\veg z)
  =
  \min_{\veg z_i\in Z_\sigma}
  \abs{\veg z_i-\veg z}
  \(
  \prod_{\veg z_i\in Z_\sigma}
  \abs{\veg z_i-\veg z}
  \)^{2/\abs{Z_\sigma}}
  \frac{1}{\abs{\veg z-\veg c_t}^{4}}
  \label{eq:tree_mimicry_selection}
\end{equation}
and for $r=1$
\begin{equation}
  \rho_1(\veg z)=1/\abs{\veg z-\veg c_t}\, .
\end{equation}
In mimicry pivoting \cite{tetzner_IncompleteAdaptiveCross_2026}, the $r$th column index is then given by
\begin{equation}
  \sigma_{\mc F_{\veg e}(t)}[r] = \argmax_{j:\veg z_j \in Z} \rho_r(\veg z_j) \, , \label{eq:mimicryargmax}
\end{equation}
where $Z$ denotes the set of positions over which the maximization is performed.
In tree mimicry pivoting, the search is accelerated by using the column tree $T_J$ itself:
Starting from the clusters contributing to $\mc F_{\veg e}(t)$, the function $\rho_r$ is evaluated for $Z = \{\veg c_s : s \subseteq \mc F_{\veg e}(t)\}$ and the search is continued recursively within the maximizing cluster (see \cite{tetzner_IncompleteAdaptiveCross_2026} for details).
Once a leaf is reached, $\rho_r$ is evaluated at the positions $\veg z_\beta$ of the associated basis functions, and the maximizing function yields the desired index $\sigma_{\mc F_{\veg e}(t)}[r]$.

The maximization in \eqref{eq:tree_mimicry_selection} is modified using the edge vectors and surface normals.
The edge vectors introduce an additional weight once a leaf is reached, whereas the surface normals restrict the sets of positions over which the maximization is performed both during the tree descent and at the basis-function level.
To this end, we associate for the row-tree traversal each \ac{RWG} $\veg \fk_\beta$ with an edge vector $\veg\ell_\beta$ and a representative surface normal $\n_\beta$.
The former is the unit vector along the supporting edge, whereas the latter is obtained by normalizing the mean of $\n$ on its two adjacent triangles.
Collinear vectors are identified by choosing a fixed representative.
The components of the resulting representatives are subsequently rounded to one decimal.

The edge vectors address the configuration in \Cref{fig:orientation_failures_edges} by introducing an additional weight in the maximization once a leaf is reached.
Let $s$ be the leaf reached at step $r$ and define
\begin{equation}
  m_{\veg \ell_\beta}^{(r)}
  =
  \abs{
    \{
    i<r :
    \veg\ell_{\sigma_{\mc F_{\veg e}(t)}[i]}
    =
    \veg\ell_\beta
    \}
  }\,,
  \qquad \beta\in s\,,
  \label{eq:edge_orientation_count}
\end{equation}
as the number of previously selected pivots associated with the edge vector $\veg\ell_\beta$.
If $\{\beta\in s:m_{\veg \ell_\beta}^{(r)}=0\}\neq\emptyset$, we set
\begin{equation}
  w_{\veg \ell_\beta}^{(r)}
  =
  \begin{cases}
    1\,, & m_{\veg \ell_\beta}^{(r)}=0\,,\\
    0\,, & \text{otherwise}.
  \end{cases}
\end{equation}
Otherwise, all edge vectors occurring in $s$ have been selected at least once and we use
\begin{equation}
  w_{\veg \ell_\beta}^{(r)}
  =
  \frac{1}{1+m_{\veg \ell_\beta}^{(r)}}\,.
  \label{eq:edge_orietntation_weight}
\end{equation}
The maximization is then performed over $\rho_r(\veg z_\beta)w_{\veg \ell_\beta}^{(r)}$.
Thus, edge vectors that have not yet been selected are initially favored exclusively, while subsequently less frequently selected edge vectors receive a larger weight.

The surface normals address the configuration in \Cref{fig:orientation_failures_normals} by restricting the set of positions over which \eqref{eq:tree_mimicry_selection} is evaluated, both during the tree descent and at the basis-function level.
For each cluster $s$ and each surface normal $\hat{\veg n}$ occurring in $s$, let
\begin{equation}
  m_{\hat{\veg n}}(s)
  =
  \abs{\{\beta\in s:\n_\beta=\hat{\veg n}\}}
  \label{eq:normal_count}
\end{equation}
denote the number of associated expansion functions.
Among the surface normals occurring in $s$, we successively select up to three normals $\hat{\veg n}_i$.
For $i=1,2,3$, we select
\begin{equation}
  \hat{\veg n}_i
  =
  \arg\max_{\substack{
      \n : \hat{\veg n}\cdot\hat{\veg n}_j=0\\
      j<i}}
  m_{\hat{\veg n}}(s)\,,
  \label{eq:normal_selection}
\end{equation}
where the orthogonality condition is omitted for $i=1$ and evaluated using the representation described above.
The selection stops if no further normal satisfying the orthogonality condition exists.
If $k_s\leq3$ normals have been selected, we use them only if
\begin{equation}
  \sum_{i=1}^{k_s} m_{\hat{\veg n}_i}(s)
  \geq
  0.75\,\abs{s}\,.
  \label{eq:normal_coverage}
\end{equation}
The remaining basis functions are collected in an additional set.
If \eqref{eq:normal_coverage} is not satisfied, we set $k_s = 0$.
The value of $75\%$ is chosen heuristically and is only intended to require that the cluster is dominated by the selected normals.

In the tree mimicry pivoting, we take the surface normals into account when selecting pivots for $\mc F_{\veg e}(t)$, only if
\begin{equation}
  k_s>0
  \qquad
  \text{for all } s\subseteq\mc F_{\veg e}(t)\,,
\end{equation}
and if the set
\begin{equation}
  \bigcup_{s\subseteq\mc F_{\veg e}(t)}
  \{\hat{\veg n}_i(s):i=1,\ldots,k_s\}
\end{equation}
contains at most three mutually orthogonal normals.
Thus, we obtain up to three $Z_n\subset Z$ with $n = 1,\dots,3$ associated with the orthogonal normals.
In that case, tree mimicry pivoting is extended by  maximizing \eqref{eq:mimicryargmax} over $j:\veg z_j \in Z_n$, where $Z_n$ alternates in round-robin order.
A set $Z_n$ is removed from the round-robin sequence only after convergence is satisfied.
After all selected normals have converged separately,  \eqref{eq:mimicryargmax} is maximized over $j:\veg z_j \in Z\setminus \bigcup_{n} Z_n$ until the convergence criterion is satisfied.
We have summarized the orientation-aware extension to tree mimicry pivoting in \Cref{alg:orientation_aware_pivoting}.

\begin{algorithm}[!t]
  \caption{Orientation-aware pivoting at step $r$}
  \label{alg:orientation_aware_pivoting}
  \begin{algorithmic}[1]
    \Procedure{OrientationAwarePivot}{$t,\mc F_{\veg e}(t)$}
    \If{surface-normal information is used}
    \If{not all selected normals have converged}
    \State $\hat{\veg n}\gets$ next unconverged  normal (round-robin)
    \State descend with $\n_\beta=\hat{\veg n}$
    \Else
    \State descend with $\n_\beta\neq\hat{\veg n}_i$ for all selected $\hat{\veg n}_i$
    \EndIf
    \Else
    \State descend unrestricted
    \EndIf
    \If{$m_{\veg\ell_\beta}^{(r)}=0$ for some $\veg z_\beta \in Z_n$}
  \State $\displaystyle
  \sigma_{\mc F_{\veg e}(t)}[r]
  =
  \arg\max_{m_{\veg\ell_\beta}^{(r)}=0}
  \rho_r(\veg z_\beta)$
  \Else
  \State $\displaystyle
  \sigma_{\mc F_{\veg e}(t)}[r]
  =
  \arg\max
  \rho_r(\veg z_\beta)\,
  w_{\veg\ell_\beta}^{(r)}$
  \EndIf
  \EndProcedure
\end{algorithmic}
\end{algorithm}

We repeatedly apply \Cref{alg:orientation_aware_pivoting}, while the row pivots are selected by the unchanged maximum-absolute-value search.
If the surface-normal information is used, \eqref{eq:iaca_convergence} is required to be satisfied separately for column indices associated with each $Z_n$, $n=1,\dots,3$ and $Z\setminus \bigcup_{n} Z_n$.
Likewise, the trend check~\cite[Eq. (28)]{tetzner_IncompleteAdaptiveCross_2026} is applied separately.
The estimate of $\normF{\mat T_{t,\sigma^{(r)}_{\mc F_{\veg e}(t)}}}$, however, is shared and updated using all selected pivots.
If the surface-normal information is not used, this reduces to the standard \ac{IACA} convergence test.

\subsection{Implementation Details and Complexity}
\label{subsec:wideband_complexity}

We first partition $\mat T$ using \eqref{eq:lfadm} and \eqref{eq:hfadm} and construct the sets of directions and parent--child mappings of \Cref{sec:DirInteractionHierarchy}.
A top-down traversal then constructs $\mc F(t)$  and $\mc F_{\veg e}(t)$ according to \eqref{eq:directional_interaction_inheritance}.
The bottom-up traversal in \Cref{alg:bottomupwidebandnca} selects the pivots and assembles the bases and transfer matrices.

\begin{algorithm}[!t]
\caption{Bottom-up wideband \ac{NCA} for the row tree}
\label{alg:bottomupwidebandnca}
\begin{algorithmic}
  \Procedure{BottomUpWidebandNCA}{$T_I$}
  \For{$\ell=L,L-1,\ldots,0$}
  \For{$t\in T_I(\ell)$}
  \If{$t$ is low frequency}
  \State select $\tau_t,\sigma_{\mc F(t)}$ by the \ac{IACA}
  \If{$t$ is leaf}
  \State assemble $\mat U_t$
  \Else
  \For{$t'_i\in S(t)$}
  \State assemble $\mat\Theta_{t'_i,t}$
  \EndFor
  \EndIf
  \Else
  \For{$\veg e\in\mc E(t)$ with $\mc F_{\veg e}(t)\neq\emptyset$}
  \State select $\tau_t^{\veg e},\sigma_{\mc F_{\veg e}(t)}$ by the \ac{IACA}
  \If{$t$ is leaf or has LF children}
  \State assemble $\mat U_t^{\veg e}$ by \eqref{eq:directional_basis}
  \Else
  \For{$t'_i\in S(t)$}
  \State $\veg e'_i\gets\veg e'(\veg e,t'_i)$ by \eqref{eq:dirmapping}
  \State assemble $\mat\Theta_{t'_i,t}^{\veg e}$ by \eqref{eq:directional_transfer}
  \EndFor
  \EndIf
  \EndFor
  \EndIf
  \EndFor
  \EndFor
  \EndProcedure
\end{algorithmic}
\end{algorithm}

After completing the row- and column-tree constructions, we assemble the coupling matrices by \eqref{eq:directional_coupling} and store the blocks in $P_\mr{nonadm}$ explicitly.
Under the standard sparsity assumption, the blocks in $P_\mr{nonadm}$ contain $\mc O(N)$ entries.

The complexity follows from the established analyses of the low- and high-frequency constructions.
For low-frequency clusters, the construction coincides with the \ac{IACA}-based \ac{NCA} for standard $\mc H^2$-matrices, whose storage, setup, and \ac{MVP} costs are $\mc O(N)$ for bounded ranks~\cite{tetzner_IncompleteAdaptiveCross_2026}.
For high-frequency clusters, the number of directions satisfies
\begin{equation}
\abs{\mc E(t)}
=
\mc O\!\left((\kappa\,\diam(B_t))^2\right)
\label{eq:direction_count_complexity}
\end{equation}
in three dimensions.
Under the standard geometric regularity assumptions for a balanced cluster tree and a wavelength-resolving discretization, summing \eqref{eq:direction_count_complexity} over the hierarchy yields
\begin{equation}
\sum_{t\in T_I}\abs{\mc E(t)}
=
\mc O(N\log N)\, ,
\label{eq:total_direction_count}
\end{equation}
and analogously for $T_J$~\cite{bebendorf_WidebandNestedCross_2015,borm_Directional$mathcalH^2$matrix_2017}.
Under these assumptions and for bounded directional ranks, the storage and \ac{MVP} costs of the high-frequency representation therefore scale as $\mc O(N\log N)$.
Since the \ac{IACA} evaluates only the matrix entries required for this representation, the setup has the same asymptotic complexity for a bounded cost per \ac{EFIE} matrix entry~\cite{tetzner_IncompleteAdaptiveCross_2026}.

The modifications introduced in this work preserve these bounds.
The spherical-Fibonacci placement changes the positions of the directions but not their asymptotic number and, therefore, retains the scaling \eqref{eq:direction_count_complexity}.
The nearest-neighbor relation \eqref{eq:dirmapping} can be evaluated with constant work per parent direction using the inverse spherical-Fibonacci mapping from \cite{keinert_SphericalFibonacciMapping_2015}, and thus linear work in the number of stored directions. The orientation-aware pivoting processes at most three normals and adds a constant amount of edge-vector bookkeeping per pivot.
At the high-to-low-frequency transitions, the high-frequency parent clusters have bounded electrical diameter and therefore $\mc O(1)$ directions by \eqref{eq:direction_count_complexity}.
Since their index sets are disjoint and the directional ranks are bounded, storing the corresponding bases $\mat U_t^{\veg e}$ requires only $\mc O(N)$ additional coefficients.
In summary, these modifications affect only constants and, hence, do not change the asymptotic complexity. 

The low-frequency cluster hierarchy contributes $\mc O(N)$ to storage, setup, and a single \ac{MVP}, while the high-frequency cluster hierarchy determines the $\mc O(N\log N)$ worst-case complexity.
Consequently, the construction attains
\begin{equation}
\begin{cases}
  \mc O(N)\,,       & \text{entirely low-frequency problems},      \\
  \mc O(N\log N)\,, & \text{problems with high-frequency clusters},
\end{cases}
\label{eq:wideband_complexity}
\end{equation}
for storage, setup, and a single \ac{MVP} under the assumptions above, consistent with the scalings measured in \Cref{sec:NumericalResults}.

\section{Numerical Results}
\label{sec:NumericalResults}

In the following examples, we study the accuracy, storage requirements, setup time, and \ac{MVP} time of the proposed wideband \ac{DH2} construction.
If not stated otherwise, storage, setup time, and \ac{MVP} time are measured for the complete compressed matrix representation, including the blocks in $P_\mr{nonadm}$, coupling matrices, leaf bases, and transfer matrices.
The setup time does not include the generation of the surface mesh.
The accuracy of a compressed matrix $\widetilde{\mat X}$ is measured by the relative error
\begin{equation}
\aleph(\widetilde{\mat X})
=
\frac{
  \norm{\widetilde{\mat X}-\mat X^\mr{ref}}_2
}{
  \norm{\mat X^\mr{ref}}_2
}\,.
\label{eq:numerical_error}
\end{equation}
Unless stated otherwise, we report $\aleph_\mr{adm}\coloneqq\aleph(\widetilde{\mat T}_\mr{adm})$, evaluated on the part of the matrix associated with admissible interactions; where the full-matrix error is reported, we write $\aleph\coloneqq\aleph(\widetilde{\mat T})$.
The spectral norm in \eqref{eq:numerical_error} is approximated by power iteration.
Unless stated otherwise, the reference $\mat X^\mr{ref}$ is obtained from an $\mc H$-matrix constructed with the random-sampling \ac{ACA} from~\cite{tetzner_AdaptiveCrossApproximation_2024} using the tolerance $10^{-5}$.

All system matrices are obtained from the \ac{EFIE} discretized with \ac{RWG} basis functions.
All compared compression methods use the admissibility conditions \eqref{eq:lfadm} and \eqref{eq:hfadm}.
Thus, the $\mc H$-matrix comparison is also based on the same block partition; it differs from the \ac{DH2} constructions only in that each admissible block is compressed independently.
For the $\mc H$-matrix comparison shown in the plots, the \ac{ACA} tolerance is chosen such that $\aleph_\mr{adm}$ is approximately equal to that of the proposed \ac{DH2}.
The proposed compression uses the orientation-aware \ac{IACA}.
The K-means cluster tree is constructed as a binary tree with $n_\mr{min}=100$.
For the octree, we use $n_\mr{min}=200$ to avoid very small child clusters after subdivision.
The directional-resolution parameter is $\gamma=1$ in all examples.
The compression tolerance of the proposed \ac{DH2} is $\varepsilon=10^{-3}$ unless stated otherwise.

All computations are carried out in double precision.
The implementation is written in Julia and uses \texttt{BEAST.jl}~\cite{krcools_KrcoolsBEASTjlV2100_2026a} for the \ac{EFIE} discretization and matrix-entry evaluation, \texttt{H2Trees.jl}~\cite{danijeljukic_Djukic14H2TreesjlV041_2026} for the cluster-tree construction, and \texttt{AdaptiveCrossApproximation.jl}~\cite{tetzner_AdaptiveCrossApproximationjl_2026a} and \texttt{NestedCrossApproximation.jl}~\cite{tetzner_NestedCrossApproximationjl_2026} for the compression algorithms.
The simulation scripts required to reproduce the numerical results are published alongside this article~\cite{tetzner_WidebandDH2CompressionforEFIE_2026}.
The timings are measured on a server with two AMD EPYC 9755 processors and $6\,\mathrm{TB}$ of memory.
The \ac{MVP} timings are single-threaded.

\subsection{Cube}
\label{subsec:numerics_cube}

First, we consider a regularly meshed cube with side length \SI{1}{\meter} as a case that is particularly well matched to octree clustering.
We decrease the wavelength from approximately \SI{0.5}{\meter} to \SI{0.1}{\meter} and use meshes with edge length $h\approx\lambda/10$.
The resulting discretizations range from $N=23\,328$ to $N=394\,272$ unknowns.
For this example, we use $\eta_{\mr l}=1$ and $\eta_{\mr h}=5$.

We compare $\mc H$-matrices constructed with the standard \ac{ACA} and with the random-sampling \ac{ACA} of~\cite{tetzner_AdaptiveCrossApproximation_2024}.
For the directional constructions, we apply the standard \ac{IACA} to the proposed K-means tree and compare it with the orientation-aware \ac{IACA} on the octree and on the proposed K-means tree.

\begin{figure}
\centering
\captionsetup[sub]{skip=0pt}
\makebox[\columnwidth][r]{%
  \subfloat[]{
    \includegraphics{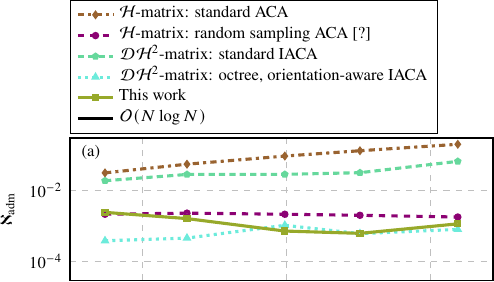}
    \label{fig:cube_relerror}
  }
}\\[-0.25cm]
\makebox[\columnwidth][r]{%
  \subfloat[]{
    \includegraphics{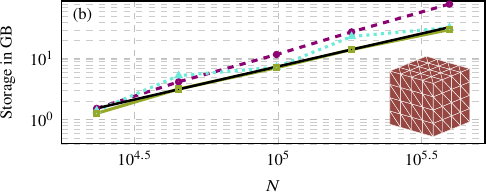}
    \label{fig:cube_storage}
  }
}\\[-8.9pt]
\caption{Cube: $\aleph_\mr{adm}$ of the considered approaches and storage requirement of the accurate approaches as functions of the number of unknowns with $h=\lambda/10$.}
\label{fig:cube_scaling}
\end{figure}

\Cref{fig:cube_relerror} shows that the standard \ac{ACA} and the standard \ac{IACA} converge prematurely and fail to reach the prescribed accuracy.
The random-sampling \ac{ACA} and both orientation-aware directional constructions yield $\aleph_\mr{adm}$ of the order of the prescribed tolerance.
The random-sampling \ac{ACA} therefore also avoids premature convergence, but it remains an $\mc H$-matrix technique and does not provide the nested, entry-efficient representation of the \ac{IACA}.
The orientation-aware \ac{IACA} restores the accuracy while retaining the \ac{DH2} construction.

Only the approaches without premature convergence are included in the storage comparison in \Cref{fig:cube_storage}.
Both \acp{DH2} require substantially less storage than the $\mc H$-matrix, while the octree and the proposed K-means tree yield similar results.
The larger octree values at the second and fourth discretizations occur because the leaf-size criterion does not yet permit an additional octree level, whereas the binary K-means tree refines more gradually.
For the largest example, the proposed construction requires \SI{30.2}{\giga\byte}, compared with \SI{32.5}{\giga\byte} for the octree construction and \SI{79.5}{\giga\byte} for the $\mc H$-matrix.
Thus, even on a geometry that is particularly favorable for an octree, the proposed geometry-adaptive construction is slightly more storage efficient, while both \acp{DH2} retain a clear advantage over the $\mc H$-matrix.

\subsection{Sphere}
\label{subsec:numerics_sphere}

Next, we consider an icosphere with radius \SI{1}{\meter}.
All sphere meshes use an average edge length of approximately $h=\lambda/10$.
We first investigate the relation between the prescribed \ac{IACA} tolerance and the resulting $\aleph_\mr{adm}$ for a discretization with approximately $10^5$ unknowns using a fully assembled matrix as reference.
The tolerance is varied from $\varepsilon=10^{-2}$ to $\varepsilon=10^{-10}$.
For this tolerance sweep, we use $\eta_{\mr l}=1$ and compare $\eta_{\mr h}=1$ and $\eta_{\mr h}=5$ for both the octree and the proposed K-means construction.
Increasing $\eta_{\mr h}$ relaxes the high-frequency admissibility condition, so that interactions become admissible earlier during the block-tree subdivision, which can substantially improve the compression efficiency.
At the same time, the resulting interactions are less strongly separated and can therefore be more difficult to approximate.

\begin{figure}
\centering
\includegraphics{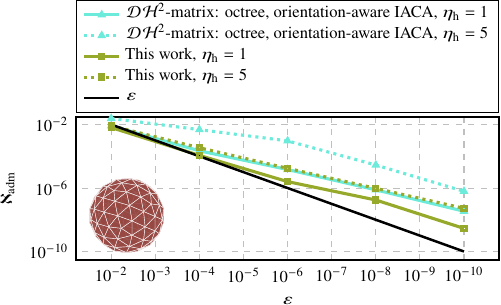}
\caption{Icosphere: $\aleph_\mr{adm}$ as a function of the \ac{IACA} tolerance for approximately $10^5$ unknowns.}
\label{fig:sphere_errorsweep}
\end{figure}

As shown in \Cref{fig:sphere_errorsweep}, $\aleph_\mr{adm}$ follows the prescribed tolerance down to approximately $\varepsilon=10^{-5}$.
For smaller tolerances, the error continues to decrease, but the gap to the prescribed tolerance increases.
The results for $\eta_{\mr h}=5$ are slightly less accurate than those for $\eta_{\mr h}=1$.
This is expected, since $\eta_{\mr h}=5$ introduces additional practically motivated heuristics in the directional construction.

We next increase the electrical size of the sphere.
For this scaling experiment, we use $\eta_{\mr l}=1$ and $\eta_{\mr h}=5$.
We compare the proposed K-means construction with the orientation-aware directional construction on the octree and with an $\mc H$-matrix.
For the $\mc H$-matrix comparison, the \ac{ACA} tolerance is adjusted such that $\aleph_\mr{adm}$ is comparable to that of the proposed construction.

\begin{figure}
\centering
\captionsetup[sub]{skip=0pt}
\makebox[\columnwidth][r]{%
  \subfloat[]{
    \includegraphics{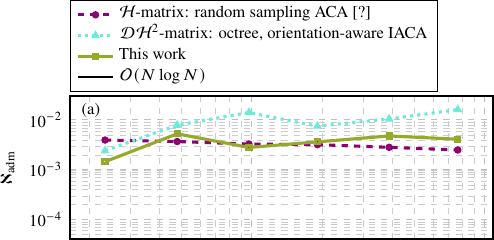}
    \label{fig:sphere_relerror}
  }
}\\[-0.25cm]
\makebox[\columnwidth][r]{%
  \subfloat[]{
    \includegraphics{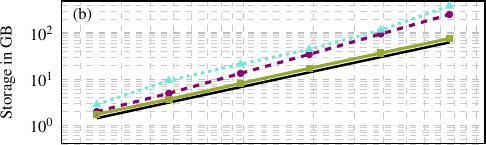}
    \label{fig:sphere_storage}
  }
}\\[-0.25cm]
\makebox[\columnwidth][r]{%
  \subfloat[]{
    \includegraphics{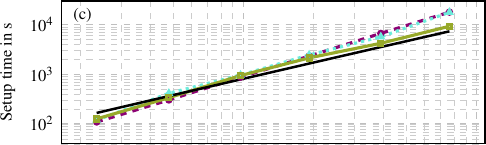}
    \label{fig:sphere_time}
  }
}\\[-0.25cm]
\makebox[\columnwidth][r]{%
  \subfloat[]{
    \includegraphics{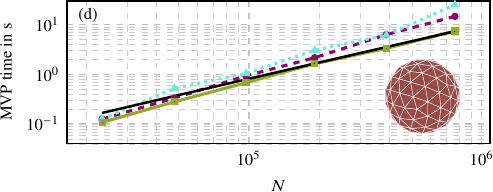}
    \label{fig:sphere_mvptime}
  }
}\\[-8.9pt]
\caption{Icosphere: $\aleph_\mr{adm}$, storage requirement, setup time, and \ac{MVP} time as functions of the number of unknowns for $h=\lambda/10$.}
\label{fig:sphere_scaling}
\end{figure}

The relative error $\aleph_\mr{adm}$ in \Cref{fig:sphere_relerror} remains of the order of the prescribed tolerance for all three constructions.
In contrast to the cube, \Cref{fig:sphere_storage} shows a clear storage advantage of the proposed K-means construction over the octree.
The difference increases with the problem size because the spherical surface is poorly matched to the box-based octree subdivision, whereas the K-means clusters follow the distribution of the unknowns on the surface more closely.
The proposed construction exhibits the expected $\mc O(N\log N)$ storage behavior.
For the largest example, with $N=768\,000$ unknowns, it requires \SI{75.9}{\giga\byte}, compared with \SI{386.1}{\giga\byte} for the octree construction and \SI{252.3}{\giga\byte} for the $\mc H$-matrix.
The timing results in \Cref{fig:sphere_time,fig:sphere_mvptime} show the same trend.
For the largest example, the proposed construction reduces the setup time to \SI{9200.8}{\second}, compared with \SI{17614.2}{\second} for the octree construction and \SI{17722.7}{\second} for the $\mc H$-matrix, and the \ac{MVP} time to \SI{7.38}{\second}, compared with \SI{24.91}{\second} and \SI{14.55}{\second}, respectively.

Finally, we validate the scattering result independently of the matrix reference by comparing the bistatic \ac{RCS} with the analytical Mie-series solution.
For this experiment, the wavelength is $\lambda=\SI{0.3}{\meter}$ and the discretization comprises $N=156\,804$ unknowns.
We use $\eta_{\mr l}=1$ and $\eta_{\mr h}=5$.
The incident plane wave propagates in the positive $z$-direction and is polarized in the $x$-direction.
The bistatic \ac{RCS} is evaluated in the $xz$-plane.
The \ac{EFIE} system is solved with \ac{GMRES}~\cite{saad_GMRESGeneralizedMinimal_1986} using a Calderón preconditioner~\cite{andriulli_MultiplicativeCalderonPreconditioner_2008}.
A relative residual of $10^{-4}$ is reached after 13 iterations.

\begin{figure}
\centering
\includegraphics{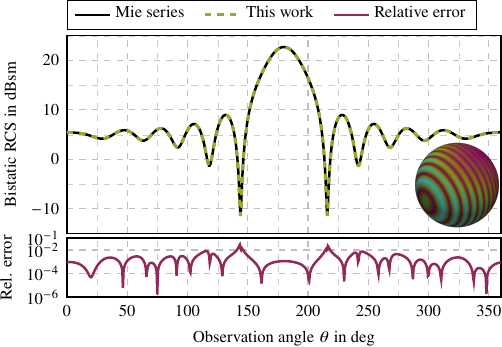}
\caption{Icosphere: bistatic \ac{RCS} in the $xz$-plane obtained with the proposed construction and the analytical Mie-series solution.}
\label{fig:sphere_mie}
\end{figure}

The numerical and analytical results in \Cref{fig:sphere_mie} agree closely over the complete observation range.
The pronounced maximum around $\theta=180^\circ$ corresponds to backscattering.
The largest relative deviations occur at the deep minima of the \ac{RCS}, where small absolute differences lead to increased relative errors.

\subsection{Rafale}
\label{subsec:numerics_rafale}

\begin{figure}
\centering
\includegraphics{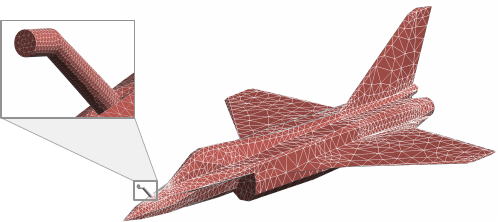}
\caption{Rafale geometry with a local mesh refinement at the nozzle.}
\label{fig:rafale_geometry}
\end{figure}

Finally, we consider the Rafale geometry shown in \Cref{fig:rafale_geometry}.
This example combines a geometry that is poorly matched to octree clustering with a strongly nonuniform discretization: 
The main surface is discretized with an average edge length of approximately $h=\lambda/10$, whereas the nozzle is discretized ten times finer with $h\approx\lambda/100$.
The resulting discretizations range from $N=35\,766$ to $N=1\,854\,723$ unknowns.

The local refinement leads to small clusters and interactions that remain in the low-frequency regime although the full geometry is electrically large.
For the Rafale scaling experiment, we use $\eta_{\mr l}=1$ and $\eta_{\mr h}=2$.
We compare the proposed K-means construction with the orientation-aware directional construction on the octree and with an $\mc H$-matrix.
For the $\mc H$-matrix comparison, the \ac{ACA} tolerance is adjusted such that $\aleph_\mr{adm}$ is comparable to that of the proposed construction.

\begin{figure}
\centering
\captionsetup[sub]{skip=0pt}
\makebox[\columnwidth][r]{%
  \subfloat[]{
    \includegraphics{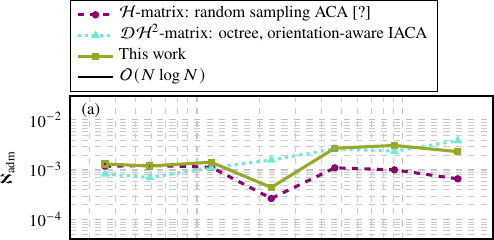}
    \label{fig:rafale_relerror}
  }
}\\[-0.25cm]
\makebox[\columnwidth][r]{%
  \subfloat[]{
    \includegraphics{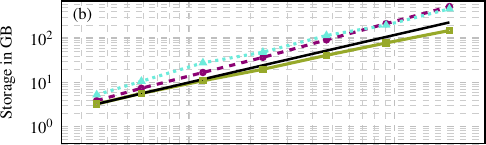}
    \label{fig:rafale_storage}
  }
}\\[-0.25cm]
\makebox[\columnwidth][r]{%
  \subfloat[]{
    \includegraphics{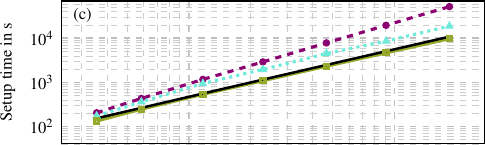}
    \label{fig:rafale_time}
  }
}\\[-0.25cm]
\makebox[\columnwidth][r]{%
  \subfloat[]{
    \includegraphics{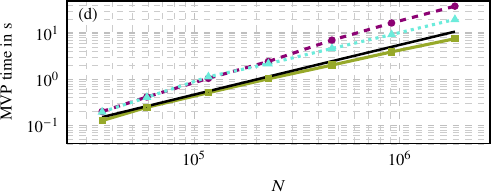}
    \label{fig:rafale_mvptime}
  }
}\\[-8.9pt]
\caption{Rafale: $\aleph_\mr{adm}$, storage requirement, setup time, and \ac{MVP} time as functions of the number of unknowns with $h_{\mr{max}}=\lambda/10$.}
\label{fig:rafale_scaling}
\end{figure}

The relative error $\aleph_\mr{adm}$ in \Cref{fig:rafale_relerror} remains of the order of the prescribed tolerance for all three constructions.
The storage results in \Cref{fig:rafale_storage} show a clear advantage of the proposed K-means construction over both the octree and the $\mc H$-matrix.
For the largest example, the proposed construction requires \SI{147.9}{\giga\byte}, compared with \SI{460.7}{\giga\byte} for the octree construction and \SI{517.8}{\giga\byte} for the $\mc H$-matrix.
The setup and \ac{MVP} times in \Cref{fig:rafale_time,fig:rafale_mvptime} show corresponding improvements.
For the largest example, the proposed construction requires \SI{9631.0}{\second} for setup and \SI{7.53}{\second} for one \ac{MVP}, compared with \SI{18605.2}{\second} and \SI{19.70}{\second} for the octree construction and \SI{51549.7}{\second} and \SI{38.20}{\second} for the $\mc H$-matrix.
Thus, the proposed construction retains its storage and timing advantages for a geometry that is poorly matched to octree clustering, a strongly nonuniform mesh, and interactions from both frequency regimes.

To illustrate the wideband transition on a fixed mesh, we finally keep the Rafale discretization with $N=742\,845$ unknowns and vary only the wavelength.
For this experiment, we use $\eta_{\mr l}=\eta_{\mr h}=1$.
In the other experiments, choosing $\eta_{\mr h}>\eta_{\mr l}$ relaxes the more restrictive high-frequency admissibility condition and can improve the compression by making high-frequency interactions admissible earlier.
Here, using the same value in both criteria deliberately removes this difference and avoids an artificial jump between the two regimes, so that the effect of the directional subdivision in the high-frequency regime can be observed more clearly.

\begin{figure}
\centering
\newcommand{\pendingaircraftpanel}[1]{%
  \fbox{%
    \parbox[c][2.8cm][c]{0.18\columnwidth}{%
      \centering
      \footnotesize
      Figure pending\\[0.4em]
      #1%
    }%
  }%
}
\IfFileExists{figures/rafale10h.png}
{\includegraphics[width=0.30\columnwidth]{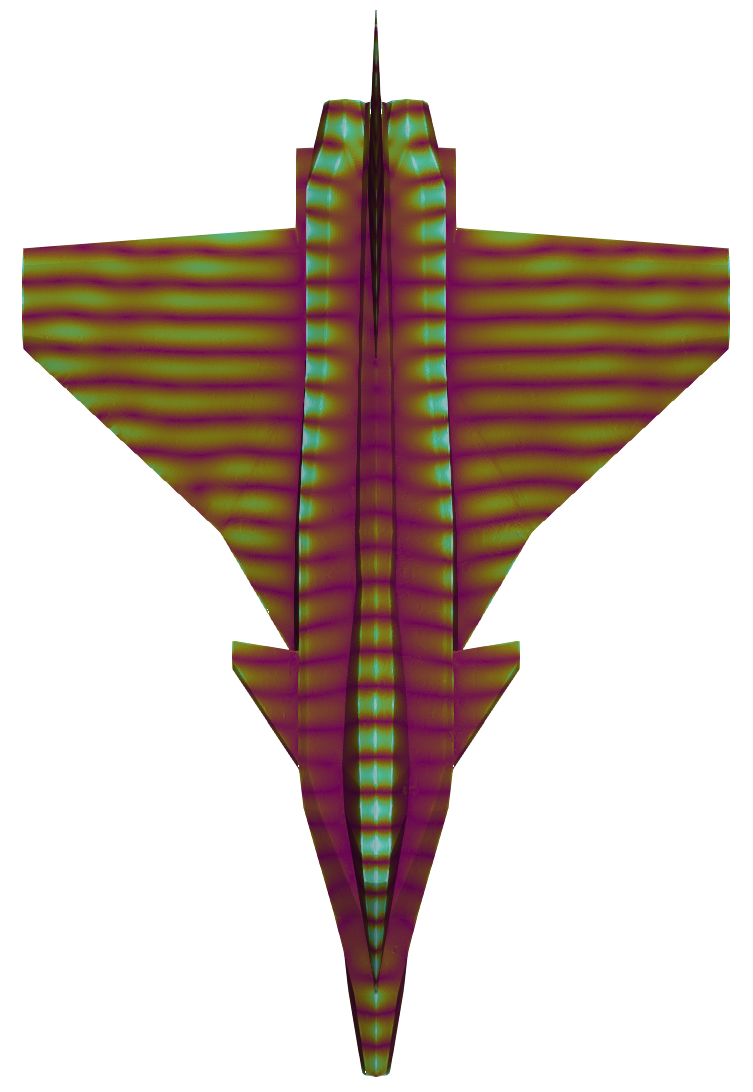}}
{\pendingaircraftpanel{Aircraft length: $16\lambda$}}%
\hfill
\IfFileExists{figures/rafale40h.png}
{\includegraphics[width=0.30\columnwidth]{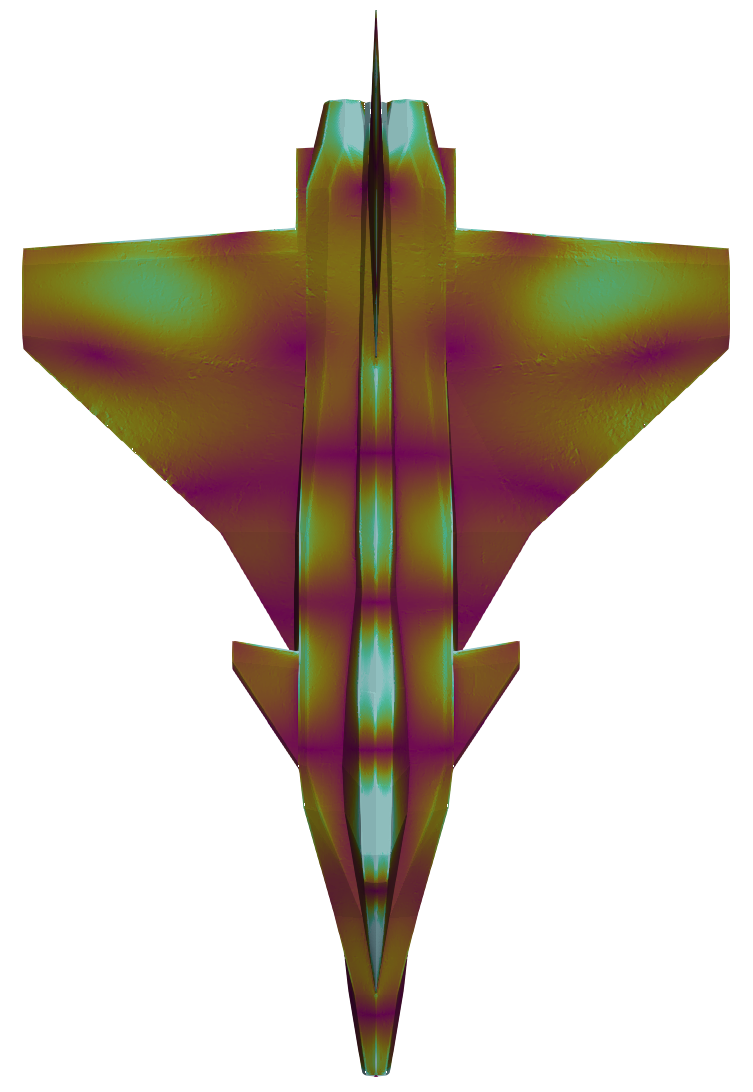}}
{\pendingaircraftpanel{Aircraft length: $4\lambda$}}%
\hfill
\IfFileExists{figures/rafale160h.png}
{\includegraphics[width=0.30\columnwidth]{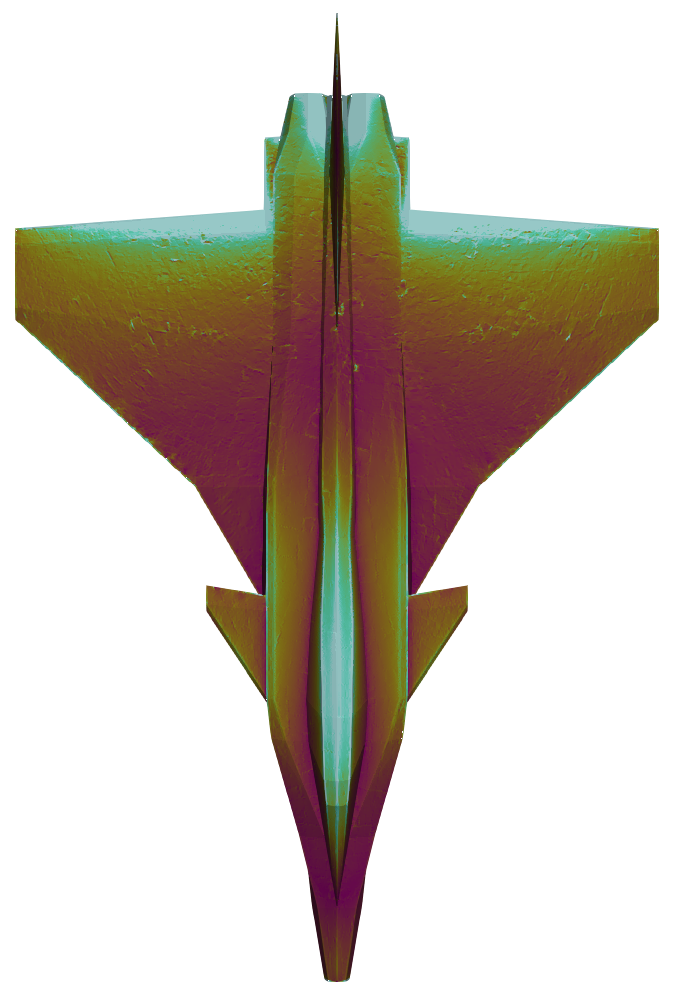}}
{\pendingaircraftpanel{Aircraft length: $\lambda$}}%
\\[0.4em]
\includegraphics{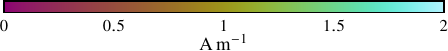}
\caption{Imaginary part of the surface current density $\veg j$ on an open aircraft model excited by a plane wave. From left to right, the aircraft lengths are approximately $16\lambda$, $4\lambda$, and $\lambda$.}
\label{fig:plane_wave_rafale}
\end{figure}

\begin{table}
\caption{Wideband behavior of the \ac{DH2} representation and iterative solution of an open aircraft model with $N=742\,845$ unknowns. The \ac{GMRES} tolerance is $10^{-6}$.}
\label{tab:rafale_wideband_comparison}
\centering
\footnotesize
\setlength{\tabcolsep}{1.5pt}
\renewcommand{\arraystretch}{1.18}
\begin{tabular*}{\columnwidth}{
    @{\extracolsep{\fill}}
    l
    *{3}{c}
    @{}
  }
  \toprule
  Aircraft length
  & $16\lambda$
  & $4\lambda$
  & $\lambda$ \\
  \cmidrule(lr){2-2}
  \cmidrule(lr){3-3}
  \cmidrule(lr){4-4}
  Low-frequency clusters
  & 5599 & 20478 & 21462 \\
  High-frequency clusters
  & 15625 & 973 & 46 \\
  Maximum number of directions
  & 37 & 13 & 3 \\
  Storage in \si{\gibi\byte}
  & 119.3 & 56.67 & 51.90 \\
  $\aleph_{\mr{adm}}$
  & $2.18\mr e{-3}$ & $2.32\mr e{-3}$ & $2.18\mr e{-3}$ \\
  $\aleph$
  & $2.08\mr e{-4}$ & $5.06\mr e{-5}$ & $4.76\mr e{-5}$ \\
  \ac{GMRES} iterations
  & 776 & 167 & 90 \\
  \bottomrule
\end{tabular*}
\end{table}

\Cref{fig:plane_wave_rafale} shows the corresponding current distributions.
As the electrical size decreases from $16\lambda$ to $\lambda$, the number of high-frequency clusters drops from 15625 to 46, while the number of low-frequency clusters increases from 5599 to 21462.
The maximum number of directions decreases from 37 to 3.
Consequently, the storage is reduced from \SI{119.3}{\gibi\byte} to \SI{51.90}{\gibi\byte}, as shown in \Cref{tab:rafale_wideband_comparison}.
The relative error $\aleph_\mr{adm}$ remains close to the prescribed tolerance, and $\aleph$ is smaller in all three cases.
The number of \ac{GMRES} iterations decreases with the electrical size, which is consistent with the less oscillatory current distributions in \Cref{fig:plane_wave_rafale}.
The experiment demonstrates that the construction changes continuously between predominantly directional high-frequency compression and predominantly non-directional low-frequency compression on the same mesh.

\section{Conclusion}

The results show that the cluster and directional hierarchies need not share a prescribed level structure.
Decoupling them allows the construction to accommodate geometries and discretizations that are poorly matched to octree clustering without changing the established complexity bounds. 
The resulting method extends the \ac{IACA} within the \ac{NCA} from electrically small problems across the low- and high-frequency regimes.


\ifCLASSOPTIONcaptionsoff
\newpage
\fi



\printbibliography

\end{document}